\documentclass{article}
\usepackage{amssymb}

\usepackage{graphicx}
\usepackage{amsmath}

\newtheorem{theorem}{Theorem}

\newtheorem{corollary}[theorem]{Corollary}

\newtheorem{definition}[theorem]{Definition}

\newtheorem{lemma}[theorem]{Lemma}

\newtheorem{proposition}[theorem]{Proposition}

\input{tcilatex}

\begin{document}

\author{m\'{e}moire sous la direction de Jean-Louis Krivine \and par Denis Bonnay}
\title{Le contenu computationnel des preuves :\\
No-counterexample interpretation\\
et sp\'{e}cification des th\'{e}or\`{e}mes de l'arithm\'{e}tique}
\date{septembre 2002\\
D.E.A. logique et fondements de l'informatique\\
Universit\'{e} Paris VII\newpage }
\maketitle
\tableofcontents

\newpage

Que nous apporte la preuve d'un th\'{e}or\`{e}me, en plus du simple fait de
savoir que le th\'{e}or\`{e}me est une cons\'{e}quence de certains axiomes ?
Il s'agit de voir quel contenu purement calculatoire peut \^{e}tre extrait
de preuves qui, lorsqu'elles utilisent toutes les ressources de la logique
classique ne vont pas forc\'{e}ment exhiber pour nous les objets dont elles
parlent. Si l'on applique ce programme \`{a} l'arithm\'{e}tique, on aimerait
que les preuves nous fournissent des fonctions calculables. Mais l'interpr%
\'{e}tation na\"{i}ve selon laquelle la preuve d'un \'{e}nonc\'{e} de la
forme $\forall _{1}\exists y_{1}...\forall x_{k}\exists
y_{k}A(x_{1}y_{1}...x_{k}y_{k})$ devrait nous fournir des fonctions
calculables $\psi _{1}(x_{1})...\psi _{k}(x_{1}...x_{k})$ telles que $%
A(x_{1}\psi _{1}(x_{1})...x_{k}\psi _{k}(x_{1}...x_{k}))$ soit vrai, ne vaut
pas, comme il est facile de le montrer. Consid\'{e}rons le pr\'{e}dicat r%
\'{e}cursif primitif $Halt(x_{1},x_{2},y)$ qui est vrai si et seulement si
la machine de Turing de num\'{e}ro $x_{1}$ avec l'entr\'{e}e $x_{2}$\ s'arr%
\^{e}te au bout de $y$ pas. On a bien :

$\vdash _{AP}\forall x_{1}\forall x_{2}\exists y(Halt(x_{1},x_{2},y)\vee
\forall z\sim Halt(x_{1},x_{2},z))$

ou encore

$\vdash _{AP}\forall x_{1}\forall x_{2}\exists y\forall
z(Halt(x_{1},x_{2},y)\vee \sim Halt(x_{1},x_{2},z))$

Mais il n'y a pas de fonction r\'{e}cursive $\psi $ telle que $%
Halt(x_{1},x_{2},\psi (x_{1},x_{2})\vee \forall z\sim Halt(x_{1},x_{2},z).$
L'interpr\'{e}tation na\"{i}ve est donc fausse au moins \`{a} partir de $\Pi
_{3}$ et donc de $\Sigma _{2}.$

Nous allons pr\'{e}senter la solution de Kreisel \`{a} ce probl\`{e}me, qui
est connue sous le nom de No-counterexample interpretation. L'id\'{e}e de la
d\'{e}monstration de Kreisel est que les preuves de non-contradiction de
l'arithm\'{e}tique, dont le but, dans l'esprit du programme Hilbertien, \'{e}%
tait de montrer que le d\'{e}tour par des moyens de preuves non-finitistes
pour prouver des \'{e}nonc\'{e}s purement num\'{e}riques \'{e}tait
acceptable, peuvent en fait \^{e}tre \'{e}galement utilis\'{e}es pour donner
un sens finitiste aux \'{e}nonc\'{e}s complexes. Nous commencerons donc par
donner la preuve de consistance d'Ackermann sur laquelle s'appuie toute la d%
\'{e}monstration de Kreisel \cite{kreis1} et \cite{kreis2}.

Nous pr\'{e}senterons enfin un des prolongements contemporains de ce
programme via la r\'{e}alisabilit\'{e} classique. Le programme d'extraction
du contenu computationnel des preuves est repris, mais les m\'{e}thodes
changent. On ne cherche plus \`{a} se ramener \`{a} des fonctions r\'{e}%
cursives, comme le faisait Kreisel ou comme le font les interpr\'{e}tations
fonctionnelles\footnote{%
Dans ces interpr\'{e}tations, c'est d'ailleurs d'abord un syst\`{e}me
intuitionniste qui est interpr\'{e}t\'{e} dans une th\'{e}orie
fonctionnelle; on cherche ensuite \`{a} \'{e}tendre l'interpr\'{e}tation
au-del\`{a} de l'arithm\'{e}tique de Peano en ajoutant des principes de r%
\'{e}currence dans la th\'{e}orie fonctionnelle (par exemple le principe de
bar-r\'{e}cursion pour $AP^{2}+AC$). Ensuite l'interpr\'{e}tation proprement
dite repose sur le pr\'{e}alable de la traduction par la double n\'{e}gation
du syst\`{e}me classique dans le syst\`{e}me intuitionniste.} dans l'esprit
de l'interpr\'{e}tation \textit{Dialectica} de G\"{o}del (voir \cite{godel}
et \cite{avfef} pour une exposition des prolongements). Le paradigme de
l'interpr\'{e}tation est celui des programmes : l'id\'{e}e derri\`{e}re les
travaux de Krivine \cite{kr2} et \cite{kr4} est d'interpr\'{e}ter
directement les axiomes comme des instructions de programmation, au sens o%
\`{u} par exemple l'extension de l'isomorphisme de Curry-Howard repose sur
l'interpr\'{e}tation de l'axiome de Pierce comme un \textit{exception-handler%
}. On pr\'{e}sente ici les r\'{e}sultat de sp\'{e}cification pour des th\'{e}%
or\`{e}mes de l'arithm\'{e}tique classique du second ordre (o\`{u} l'axiome
de compr\'{e}hension d\'{e}coule de la r\`{e}gle pour l'\'{e}limination des
quantificateurs universels du second ordre), en exhibant le lien avec
l'interpr\'{e}tation de Kreisel et en insistant sur la question de la
modularit\'{e} de l'interpr\'{e}tation.

\section{La preuve de consistance d'Ackermann}

\subsection{L'epsilon-calcul et l'arithm\'{e}tique epsilon}

L'$\epsilon $-calcul est un langage sans quantificateur dans lequel les
quantificateurs sont remplac\'{e}s par des $\epsilon $-termes que l'on peut
voir comme des fonctions de choix. Si $A(x,y_{1}...y_{n})$ est une formule
avec $x,y_{1}...y_{n}$ comme variables libres, $\epsilon
_{x}A(x,y_{1}...y_{n})$ est un terme avec $y_{1}...y_{n}$ comme variable
libres, qu'il faut voir comme une fonction donnant pour des $b_{1}...b_{n}$
quelconques un t\'{e}moin pour $A(x,b_{1}...b_{n})$ s'il y en a un.\ Les $%
\epsilon $-termes sont gouvern\'{e}s par un sch\'{e}ma d'axiome, l'axiome du
transfini de Hilbert :

$A(t)\rightarrow A(\epsilon _{x}Ax)$

On peut alors d\'{e}finir les quantificateurs de la mani\`{e}re suivante :

$\exists xA(x)\equiv A(\epsilon _{x}Ax)$

$\forall xA(x)\equiv A(\epsilon _{x}\sim Ax)$

et retrouver \`{a} partir de l\`{a} les r\`{e}gles habituelles pour $\forall
,\exists $ dans le syst\`{e}me choisi.

On rappelle ici le syst\`{e}me \`{a} la Hilbert utilis\'{e} par Ackermann,
qu'on d\'{e}signera par PA$_{\epsilon }$. A quelques variantes pr\`{e}s, il
s'agit du syst\`{e}me Z$_{\mu }$ de Hilbert et Bernays\ \cite{hilb2}. Les
termes et les formules du langage sont form\'{e}s \`{a} partir de variables,
de 0, des fonctions $d,^{\prime },+,\times ,$ du signe d'\'{e}galit\'{e} et
des symboles logiques $\rightarrow $ et $\sim $. La seule r\`{e}gle
utilis\'{e}e est le \textit{modus ponens}. Il y a trois groupes de
sch\'{e}mas d'axiomes :

1) Sch\'{e}mas d'axiomes pour le niveau propositionnel

A,B d\'{e}signent des formules sans variables libres

I.1$\qquad A\rightarrow (B\rightarrow A)$

I.2$\qquad (A\rightarrow (B\rightarrow C))\rightarrow ((A\rightarrow
B)\rightarrow (A\rightarrow C))$

I.3$\qquad (\sim A\rightarrow \sim B)\rightarrow (B\rightarrow A)$

2) Sch\'{e}mas d'axiomes de l'arithm\'{e}tique

a,b,c d\'{e}signent des termes sans variables libres

II.01\qquad $a=a$

II.02$\qquad a^{\prime }=b^{\prime }\rightarrow a=b$

II.03$\qquad a\neq 0\rightarrow d(a^{\prime })=a$

II.04$\qquad a+0=a$

II.05$\qquad a+b^{\prime }=(a+b)^{\prime }$

II.06$\qquad a\times 0=0$

II.07$\qquad a\times b^{\prime }=(a\times b)+a$

II.08$\qquad a=b\rightarrow a^{\prime }=b^{\prime }$

II.09$\qquad a=b\rightarrow d(a)=d(b)$

II.10\qquad $a=b\rightarrow a+c=b+c$

II.11\qquad $a=b\rightarrow c+a=c+b$

II.12\qquad $a=b\rightarrow a\times c=b\times c$

II.13\qquad $a=b\rightarrow c\times a=c\times b$

3) Sch\'{e}mas d'axiomes pour les $\epsilon $-termes

a,b d\'{e}signent des termes sans variable libre

III.1\qquad $A(a)\rightarrow A(\epsilon _{x}Ax)$

III.2\qquad $A(a)\rightarrow \epsilon _{x}Ax\neq a^{\prime }$

III.3\qquad $\sim A(\epsilon _{x}Ax)\rightarrow \epsilon _{x}Ax=0$

III.4\qquad $a=b\rightarrow \epsilon _{x}A(x,a)=\epsilon _{x}A(x,b)$

On remarque que toutes les formules d'une d\'{e}duction sont des formules
sans variables libres.

On aurait pu utiliser \`{a} la place du groupe I un autre syst\`{e}me
d\'{e}ductif, mais l'essentiel est ici que la d\'{e}monstration puisse se
faire avec un syst\`{e}me de logique classique : les seules
propri\'{e}t\'{e}s utilis\'{e}es dans la d\'{e}monstration de la consistance
sont que ces r\`{e}gles pr\'{e}servent la v\'{e}rit\'{e} des formules sans
quantificateurs et qu'interviennent seulement des formules closes. Les
sch\'{e}mas d'axiomes du groupe II pourraient \^{e}tre remplac\'{e}es
n'importe quels autres sch\'{e}mas correspondant \`{a} des formules
universelles vraies de l'arithm\'{e}tique. En particulier, on pourrait
enrichir le langage avec d'autres fonctions r\'{e}cursives et leurs
d\'{e}finitions ou int\'{e}grer par exemple le sch\'{e}ma de r\'{e}currence
primitive.

Les axiomes de la forme III.1-III.4 sont appel\'{e}s formules critiques.Dans
le cadre de l'arithm\'{e}tique, les $\epsilon $-termes sont vus non
seulement comme des t\'{e}moins mais comme les plus petits t\'{e}moins
possibles. Le sch\'{e}ma d'induction compl\`{e}te est d\'{e}rivable dans
l'arithm\'{e}tique-epsilon, en utilisant essentiellement III.1 et III.2. A
partir de $A(0)$ et $A(t)\rightarrow A(t^{\prime }),$ l'id\'{e}e est de
montrer $A(\epsilon _{x}\sim Ax)$, en utilisant le fait que $\epsilon
_{x}\sim Ax$ ne peut \^{e}tre diff\'{e}rent de 0.

\subsection{La m\'{e}thode de substitution}

Le but de l'article d'Ackermann est de fournir une preuve de la consistance
de l'arithm\'{e}tique-epsilon, montrant ainsi que le d\'{e}tour par le
transfini via les $\epsilon $-termes n'introduit pas de contradiction.
L'id\'{e}e est de remplacer les $\epsilon $-termes d'une preuve par des
entiers de mani\`{e}re \`{a} ce que tous les axiomes critiques soient vrais.
Comme les axiomes obtenus \`{a} partir des sch\'{e}mas de I et II sont
vrais, toutes les \'{e}tapes de la d\'{e}monstration deviennent alors des
formules purement num\'{e}riques vraies. Ceci montre qu'il n'y a pas de
preuve de $0=1 $, car cette formule est fausse. Le point qui nous
int\'{e}resse dans cette d\'{e}monstration est la possibilit\'{e} de borner
les valeurs susbstitu\'{e}es aux $\epsilon $-termes en fonction de
constantes de la preuve, mais avant d'en arriver l\`{a}, il faut donner la
preuve de non-contradiction.

\begin{definition}
On d\'{e}signera par cat\'{e}gorie (\textit{Grundtyp}) d'un $\epsilon $%
-terme clos $\epsilon _{x}Ax$ le terme obtenu \`{a} partir de celui-ci en
rempla\c{c}ant tous ces sous-termes imm\'{e}diats dans lesquels $x$
n'appara\^{i}t pas libre par des variables fra\^{i}ches.
\end{definition}

Par exemple $\epsilon _{x}(0^{\prime }+\epsilon _{y}(y=0^{\prime \prime
})=\varepsilon _{z}(z^{\prime }=x)$ appartient \`{a} la cat\'{e}gorie $%
\epsilon _{x}(w=\varepsilon _{z}(z^{\prime }=x),$ tandis que $\epsilon
_{x}(\epsilon _{y}(y=0^{\prime \prime })+x=x^{\prime \prime })$ est \`{a}
lui-m\^{e}me sa propre cat\'{e}gorie. Deux $\epsilon $-termes peuvent
appartenir \`{a} la m\^{e}me cat\'{e}gorie. Le rang (\textit{Rang)} d'une
cat\'{e}gorie est le nombre d'$\epsilon $-termes enchass\'{e}es qu'elle
contient. Le rang d'un $\epsilon $-terme est le rang de sa cat\'{e}gorie.

\begin{definition}
Une substitution (\textit{Gesammtersetzung}) relativement \`{a} un ensemble
de formules closes est l'assignation \`{a} chaque cat\'{e}gorie ayant k
variables libres de l'ensemble des cat\'{e}gories des formules d'une
fonction k-aire (aux cat\'{e}gories sans variables libres sont donc
assign\'{e}s des entiers).
\end{definition}

On appelle la fonction assign\'{e}e \`{a} une cat\'{e}gorie son substituant.
La solution d'un ensemble de formules closes relativement \`{a} une
substitution pour l'ensemble des cat\'{e}gories de ses $\epsilon $-termes,
est le r\'{e}sultat obtenu en rempla\c{c}ant dans les formules les $\epsilon
$-termes par les fonctions qui sont assign\'{e}es aux cat\'{e}gories
auxquelles ils appartiennent. On dira que le substituant d'une cat\'{e}gorie
est nul s'il s'agit de 0 ou d'une fonction constante \'{e}gale \`{a} 0.

On dira qu'une substitution pour un ensemble de cat\'{e}gories a la
propri\'{e}t\'{e} P dans le cas suivant. Si \`{a} $\epsilon
_{x}A(x,w_{1}...w_{n})$ une des cat\'{e}gories est assign\'{e}e la fonction $%
f(w_{1}...w_{n})$, alors, pour des entiers quelconques $b_{1}...b_{n},$ soit
$f(b_{1}...b_{n})=0,$ soit $A(f(b_{1}...b_{n}),b_{1}...b_{n})$ est le cas et
$f(b_{1}...b_{n})$ est le plus petit entier m tel que $A(m,b_{1}...b_{n})$.

Comme les $\epsilon $-termes d'un axiome de la forme III.4 appartiennent
\`{a} la m\^{e}me cat\'{e}gorie, une substitution rend forc\'{e}ment tous
les axiomes de cette forme vraie. De plus, si une substitution a la
propri\'{e}t\'{e} P, les axiomes de la forme III.2 et III.3 sont
\'{e}galement rendus vrais.

Relativement \`{a} une preuve vue comme ensemble de formules, on se donne
une \'{e}num\'{e}ration de ces cat\'{e}gories telle que si un $\epsilon $%
-terme $b$\ avec $x$ comme variable libre est un sous-terme d'un autre $%
\epsilon $-terme $\epsilon _{x}Ax,$ la cat\'{e}gorie de $b$\ pr\'{e}c\`{e}de
la cat\'{e}gorie de $\epsilon _{x}Ax$ (*). On d\'{e}finit alors une suite de
substitutions $S_{0}...S_{n}...$ de la mani\`{e}re suivante :

1) $S_{0}$ assigne \`{a} toutes les cat\'{e}gories des substituants nuls.

2) Si on est arriv\'{e} \`{a} la subsitution $S_{n}$ qui a la
propri\'{e}t\'{e} P, on regarde le premier axiome de la forme III.1 qui est
rendu faux par $S_{n}$. Disons qu'il s'agit de $A(a,b)\rightarrow A(\epsilon
_{x}A(x,b),b)$.\ On d\'{e}finit alors $S_{n+1}$ \ ainsi. Si la formule est
fausse, cela signifie que $A(a,b)$ est vrai et $A(\epsilon _{x}A(x,b),b)$
faux. On modifie alors la valeur du subsituant $f(x)$ de la cat\'{e}gorie de
$\epsilon _{x}A(x,b)$ pour la valeur m, o\`{u} m est la valeur donn\'{e}e
par $S_{n}$ au terme $b$, en la rempla\c{c}ant par le plus petit $n$ tel que
$A(n,m).$ Toutes les cat\'{e}gories qui suivent la cat\'{e}gorie de $%
\epsilon _{x}A(x,b)$ re\c{c}oivent des subsituants nuls. Les autres
substituants sont inchang\'{e}s. On v\'{e}rifie que $S_{n+1}$ a encore la
propri\'{e}t\'{e} P. Soit $\epsilon _{x}B(x,w_{1}...w_{n})$ une
cat\'{e}gorie. Si elle suit la cat\'{e}gorie dont on a modifi\'{e} le
substituant, son nouveau substituant est nul et donc la condition est
trivialement remplie. Si elle la pr\'{e}c\`{e}de, ses substituants, ainsi
que tous ceux d'autres $\epsilon $-termes de B $B(x,w_{1}...w_{n})$ \`{a}
cause de (*) sont inchang\'{e}s, donc la propri\'{e}t\'{e} est
h\'{e}rit\'{e}e de $S_{n}$. Quant \`{a} $\epsilon _{x}A(x,y)$ elle-m\^{e}me,
la condition est v\'{e}rifi\'{e}e lorsque $y$ re\c{c}oit $b$, par
construction et par (*), et pour les autres elle est v\'{e}rifi\'{e}e comme
dans le cas pr\'{e}c\'{e}dent.

Lorsqu'on applique ces substitutions aux formules d'une preuve, parmi les
axiomes utilis\'{e}s, seuls ceux de la forme III.1 sont susceptibles
d'\^{e}tre faux. Si la suite de substitutions s'arr\^{e}te pour un certain $%
n,$ c'est donc que tous les axiomes du r\'{e}sultat sont vrais, et toutes
les formules du r\'{e}sultat \'{e}galement. Il suffit donc de montrer que la
suite de substitutions termine. Intuitivement, il semble bien qu'on corrige
peu \`{a} peu la substitution initiale. N\'{e}anmoins, il serait erron\'{e}
de penser que chaque substitution corrige une fois pour toutes un axiome de
sorte que la suite des substitutions finisse par \'{e}puiser la suite des
axiomes. Par exemple, en modifiant le substituant de $\epsilon _{x}A(x,y)$
pour rendre vrai $A(a,b)\rightarrow A(\epsilon _{x}A(x,b),b),$ il se peut
que je rende un autre axiome $B(c,d)\rightarrow B(\epsilon _{x}B(x,d),d)$
qui \'{e}tait vrai en modifiant la valeur de $d$.

\subsection{La preuve de terminaison}

\begin{definition}
Si une preuve a $m$ cat\'{e}gories et qu'une substitution $S$\ annule tous
les substituants \`{a} partir de la $n$-i\`{e}me cat\'{e}gorie, on attribue
\`{a} $S$ le nombre caract\'{e}ristique $m-s.$
\end{definition}

\begin{definition}
On ordonne les $\epsilon $-termes clos de la preuve, de mani\`{e}re \`{a} ce
qu'un sous-terme pr\'{e}c\`{e}de le terme dans lequel il appara\^{i}t. Soit $%
a_{0}...a_{k}$ une telle suite de termes, on d\'{e}finit l'ordre (\textit{%
Reduktionsgrad}) d'une substitution S par $%
o(S)=2^{k}f(o)+2^{k-1}f(1)+...+2^{0}f(k)$ o\`{u} $f(i)=1$ si $S$ donne \`{a}
$a_{i}$ la valeur 0 et $f(i)=0$ sinon.
\end{definition}

\begin{definition}
Le degr\'{e} d'une substitution S, not\'{e} $d(S)$,\ est son ordre
relativement, non plus aux formules de la preuve, mais \`{a} l'ensemble de
formules $A(0,b)...A(k,b)$ o\`{u} $A(a,b)\rightarrow A(\epsilon
_{x}A(x,b),b) $ est le premier axiome rendu faux par S et o\`{u} $k$ est la
valeur que S donne \`{a} $a$. Si S est la substitution finale, on fixe $%
d(S)=0.$
\end{definition}

\begin{definition}
L'indice, not\'{e} $i(S)$, d'une substitution S est alors la paire (o,d)
o\`{u} o est l'ordre et d le degr\'{e} de S. On ordonne les indices selon
l'ordre lexicographique.
\end{definition}

\begin{definition}
On dira qu'une substitution T est (strictement) progressive sur une
substitution S si chaque fois qu'un substituant de S prend une valeur non
nulle, T lui donne la m\^{e}me valeur (et S n'est pas progressive sur T).
\end{definition}

\begin{theorem}
Soient S et T deux substitutions fournissant des substituants pour toutes
les cat\'{e}gories d'un ensemble de formules closes, si T est progressive
sur S, soit 1) $o(T)<o(S)$, soit 2) les $\epsilon $-termes des formules re\c{%
c}oivent tous les m\^{e}mes valeurs par S et T.
\end{theorem}

Sous l'hypoth\`{e}se de progressivit\'{e}, on montre non 2) implique 1).
Soit $a_{i}$ le premier $\epsilon $-terme dans l'\'{e}num\'{e}ration sur
lequel S et T diff\`{e}rent; comme T est progressive sur S et que T et S
s'accordent sur tous les $\epsilon $-termes pr\'{e}c\'{e}dents, cela veut
dire que S annule $a_{i}.$ Par d\'{e}finition de $o,$ $o(T)<o(S).$

\begin{theorem}
Si S$_{j}$ est progressive sur S$_{i},$ alors soit 1) $i(S_{j})<i(S_{i}),$
soit 2) S$_{j+1}$ est progressive sur S$_{i+1}$ et elles sont obtenues de la
m\^{e}me mani\`{e}re \`{a} partir de S$_{j}$ et S$_{i}.$
\end{theorem}

Si $o(S_{j})<o(S_{i}),$ 1) est le cas. Sinon, par le th\'{e}or\`{e}me
pr\'{e}c\'{e}dent, les $\epsilon $-termes re\c{c}oivent les m\^{e}mes
valeurs par les deux substitutions. Par cons\'{e}quent, le premier axiome $%
A(a,b)\rightarrow A(\epsilon _{x}A(x,b),b)$ qu'elles rendent faux est le
m\^{e}me. On consid\`{e}re alors l'ensemble de formules $A(0,b)...A(k,b)$. S$%
_{j}$ est toujours progressive sur S$_{i}$ relativement \`{a} cet ensemble.
On r\'{e}applique le th\'{e}or\`{e}me pr\'{e}c\'{e}dent. Si $%
d(S_{j})<d(S_{i}),$ on a $i(S_{j})<i(S_{i});$ sinon $\epsilon _{x}A(x,b)$
recevra la m\^{e}me valeur dans $S_{j}$ et $S_{i}.$

\begin{definition}
On d\'{e}finit la notion de m--s\'{e}rie (\textit{m-Reihe)} de
substitutions. Une 1-s\'{e}rie est constitu\'{e}e d'une seule substitution.
Une m+1-s\'{e}rie est une suite $S_{i}S_{i+1}...S_{i+t}$ ($t\geq 0$) de
substitutions successives telles que les nombres caract\'{e}ristiques de $%
S_{i+1}...S_{i+t}$ sont $\leq m$ et ceux de $S_{i},S_{i+t+1}$ sont $\geq m+1$
(ou bien $S_{i+t}$ est la substitution finale).
\end{definition}

On remarque qu'une m+1-s\'{e}rie est constitu\'{e}e d'au moins une
m-s\'{e}rie.

\begin{definition}
L'indice d'une m-s\'{e}rie est un ordinal d\'{e}fini \ de la mani\`{e}re
suivante. Une 1-s\'{e}rie a l'indice $\omega o+d$ o\`{u} $(o,d)$ est le
couple d'entiers attribu\'{e} \`{a} la substitution qui la constitue. Une
m+1 s\'{e}rie est constitu\'{e}e d'une suite de m-s\'{e}ries successives
d'indices $a_{1}...a_{n},$ elle a pour indice $\omega ^{a_{1}}+...+\omega
^{a_{n}}.$
\end{definition}

\begin{theorem}
Soient $S_{1}...S_{j},T_{j+1}...T_{j+l+1}$ deux m-s\'{e}ries
cons\'{e}cutives avec $T_{j+1}$ de nombre caract\'{e}ristique $m$ et $%
a_{1}...a_{j},b_{j+1}...b_{j+l+1}$ les indices des 1-s\'{e}ries qui les
constituent. 1) $T_{j+1}$ est strictement progressive sur $S_{1}$ et 2) On
peut trouver un entier p tel que $a_{p}>b_{j+p},$ et pour $1\leq i<p,$ $%
a_{p}=b_{j+p},$ et $S_{i}$ et $T_{j+1}$ ont le m\^{e}me nombre
caract\'{e}ristique pour $1<i\leq p.$
\end{theorem}

1) Soit $g$ le nombre de cat\'{e}gories. On consid\`{e}re la suite de ces $g$
cat\'{e}gories $\epsilon _{1}...\epsilon _{g-(m+1)},\epsilon
_{g-m}...\epsilon _{g}.$ On sait d\'{e}j\`{a} que $S_{1}$ et $T_{j+1}$
annulent les substituants de $\epsilon _{g-m}...\epsilon _{g}.$ Maintenant,
comme $S_{2}...S_{j}$ ont des nombres caract\'{e}ristiques $<m$, elles ne
touchent pas aux substituants de $\epsilon _{1}...\epsilon _{g-(m+1)}.$
Comme $T_{j+1}$ est de nombre caract\'{e}ristique $m$, elle est obtenue
\`{a} partir de $S_{j}$ en rendant positive une valeur nulle de $\epsilon
_{g-(m+1)}.$ Donc $T_{j+1}$ est bien strictement progressive sur $S_{1}.$

2) Si $a_{1}>b_{j+1},$ il suffit de prendre p=1. $a_{1}<b_{j+1}$ car $%
T_{j+1} $ est bien strictement progressive sur $S_{1}.$ On se place
maintenant dans la configuration $a_{1}=b_{j+1}.$ On va supposer qu'il n'y a
pas de p tel que $a_{p}\neq b_{j+p}$ et en d\'{e}river une contradiction.

Premier cas, la premi\`{e}re m-s\'{e}rie est strictement plus longue. On
remarque que $T_{j+l+1}$ ne peut \^{e}tre la substitution finale, car par
hypoth\`{e}se, elle a le m\^{e}me indice $(o,d)$ que $S_{l+1},$ de sorte
qu'on aurait alors $d($ $S_{l+1})=0$ de sorte que $S_{l+1}$ aurait
d\'{e}j\`{a} d\^{u} \^{e}tre la substitution finale. On montre facilement
par induction que $T_{j+l+1}$ est progressive sur $S_{l+1}.$ Pour $1,$ cela
d\'{e}coule de 1). Pour $k+1$ ($k\leq l$), cela d\'{e}coule de
l'hypoth\`{e}se de r\'{e}currence et du fait que $T_{j+k}$ et $S_{k}$ ont
m\^{e}me indice, de sorte que le passage \`{a} $T_{j+k+1}$ et $S_{k+1}$
correspond \`{a} la m\^{e}me modification. Comme $T_{j+l+1}$ est progressive
sur $S_{l+1},$ on peut appliquer le th\'{e}or\`{e}me 9 et en d\'{e}duire que
$T_{j+l+2}$ est progressive sur $S_{l+2}.$ Or ceci est impossible car $%
S_{l+2} $ a un nombre caract\'{e}ristique $<m$ tandis que comme $T_{j+l+2}$
suit une m-s\'{e}rie, elle a un nombre caract\'{e}ristique $\geq m.$ On peut
donc trouver un $p$ tel que $a_{p}\neq b_{j+p}.$

Deuxi\`{e}me cas : la deuxi\`{e}me m-s\'{e}rie est au moins aussi longue que
la premi\`{e}re. $S_{j}$ et $T_{2j+1}$ ont m\^{e}me indice par
hypoth\`{e}se. Donc $T_{j+1}$ et $T_{2(j+1)}$ re\c{c}oivent le m\^{e}me
nouveau substituant positif. Mais c'est impossible car $T_{2j+1}$ a un
nombre caract\'{e}ristique $<m,$ ce qui veut dire qu'il a pr\'{e}serv\'{e}
le nouveau substituant positif de $T_{j+1}.$

On a donc d\'{e}montr\'{e} l'existence d'un premier $p$ tel que $a_{p}\neq
b_{j+p}.$ Et on doit bien avoir alors $a_{p}>b_{j+p}$ \`{a} cause du fait
d\'{e}montr\'{e} \`{a} l'instant que $T_{j+k}$ est progressive sur $S_{k}$
pour $k\leq p.$ Et $S_{i}$ et $T_{j+1}$ ont bien le m\^{e}me nombre
caract\'{e}ristique pour $1<i\leq p$, \'{e}tant donn\'{e} que les
r\'{e}sultats des substitutions sont les m\^{e}mes (par progressivit\'{e} et
\'{e}galit\'{e} des indices pour $i<p$). On appellera $S_{p}$ et $T_{j+1+p}$
les substitutions d\'{e}termin\'{e}es par les deux m-s\'{e}ries
cons\'{e}cutives.

\begin{theorem}
On se donne les m\^{e}mes conditions que dans le th\'{e}or\`{e}me
pr\'{e}c\'{e}dent ainsi qu'un $p$ tel que $1\leq p\leq m$ et les indices $%
a_{1}...a_{s},$ $b_{1}...b_{t}$ des p-s\'{e}ries $\Sigma _{1}...\Sigma _{s},$
$\Upsilon _{1}...\Upsilon _{t}$qui composent les deux m-s\'{e}ries. Il
existe alors un entier $q$ tel que $a_{q}$ et $b_{q}$ sont les indices des
deux p-s\'{e}ries qui contiennent les substitutions d\'{e}termin\'{e}es. De
plus, $a_{q}>b_{q}$ et $a_{q^{\prime }}=b_{q^{\prime }}$ pour $q^{\prime
}<q. $
\end{theorem}

La d\'{e}monstration se fait par r\'{e}currence sur $(m,p)$. Le
th\'{e}or\`{e}me pr\'{e}c\'{e}dent nous donne le cas de base et les cas $m$
quelconque et $p=1.$ On suppose maintenant que le r\'{e}sultat vaut pour un $%
m$ fix\'{e} jusqu'\`{a} $p$ exclus ($p>1$) et pour tous les $(m^{\prime
},p^{\prime })$ o\`{u} $m^{\prime }<m.$ On veut le r\'{e}sultat pour $(m,p).$
Le th\'{e}or\`{e}me pr\'{e}c\'{e}dent nous permet de d\'{e}terminer les
substitutions d\'{e}termin\'{e}es $S_{p}$ et $T_{j+1+p}.$ Pour $i<p,$ $S_{i}$
et $T_{j+1+i}$ ont les m\^{e}mes indice et nombre caract\'{e}ristique, donc
il existe un $q$ tel que $S_{p}$ appartient \`{a} $\Sigma _{q}$ et $%
T_{j+1+q} $ appartient \`{a} $\Upsilon _{q}$ tel que pour $j<q,$ $%
a_{j}=b_{j}.$ On d\'{e}compose $\Sigma _{q}$ et $\Upsilon _{q}$ en $p-1$%
-s\'{e}ries $\Gamma _{1}...\Gamma _{s^{\prime }},$ $\Delta _{1}...\Delta
_{t^{\prime }}$ d'indices $c_{1}...c_{s^{\prime }},$ $d_{1}...d_{t^{\prime }}
$. On applique l'hypoth\`{e}se de r\'{e}currence pour $(m,p-1)$. Le
d\'{e}crochage se fait n\'{e}cessairement \`{a} l'int\'{e}rieur de $\Sigma
_{q}$ et $\Upsilon _{q}$, donc on trouve un $q^{\prime }$ tel que $%
c_{q^{\prime }}>d_{q^{\prime }}$ et pour $j^{\prime }<q^{\prime },$ $%
c_{j^{\prime }}=d_{j^{\prime }}.$ Donc $a_{q}=\omega ^{c_{1}}+...+\omega
^{c_{q^{\prime }-1}}+\omega ^{c_{q^{\prime }}}+...+\omega ^{c_{s^{\prime }}}$
et $b_{q}=\omega ^{c_{1}}+...+\omega ^{c_{q^{\prime }-1}}+\omega
^{d_{q^{\prime }}}+...+\omega ^{d_{t^{\prime }}}. $ En appliquant
l'hypoth\`{e}se de r\'{e}currence avec $(p-1,p-1)$, on a que $%
c_{1}>c_{2}>...>c_{s^{\prime }}$ et$\ d_{1}>d_{2}...>d_{t^{\prime }}.$ On
peut alors conclure $a_{q}>b_{q}.$

\begin{corollary}
Deux m-s\'{e}ries cons\'{e}cutives $S_{1}...S_{j},T_{j+1}...T_{j+l+1}$ avec $%
T_{j+1}$ de nombre caract\'{e}ristique $m$ ont des indices strictement
d\'{e}croissants.
\end{corollary}

\begin{theorem}
Pour toute preuve, la suite des substitutions termine sur un r\'{e}sultat
qui rend vraie toutes les formules de la preuve.
\end{theorem}

Soit g le nombre de cat\'{e}gories de la preuve. Comme il n'y a pas de
cha\^{i}ne infinie descendante d'ordinaux, il y a un nombre fini de
g-s\'{e}ries de la forme $S_{i}...S_{i+t}$ avec $S_{i+t}$ non finale.
Apr\`{e}s cela viennent un nombre fini de g-1 s\'{e}ries et ainsi de suite
jusqu'\`{a} la substitution finale qui est la derni\`{e}re 1-s\'{e}rie. La
suite des substitutions termine donc, ce qui veut dire que tous les axiomes
IIII.1 sont rendus vrais; comme tous les autres axiomes sont par ailleurs
rendus vrais par toutes les substitutions de la suite, toutes les formules
de la preuve sont \textit{in fine} rendues vraies.

\subsection{Le bornage}

La preuve de terminaison nous a montr\'{e} qu'\'{e}tant donn\'{e}e une
preuve dans l'arithm\'{e}tique-epsilon, on pouvait trouver des substituants,
sous la forme de fonctions r\'{e}cursives (et m\^{e}me de fonctions toujours
\'{e}gales \`{a} 0 sauf un nombre fini de fois) aux cat\'{e}gories de
mani\`{e}re \`{a} ce que dans le r\'{e}sultat de la substitution toutes les
formules soient vraies. Dans le cas de la preuve d'une formule
existentielle, ce qui correspond dans le langage epsilon \`{a} une formule
de la forme $A(\epsilon _{x}Ax),$ on obtient donc une substitution qui
r\'{e}sout le terme $\epsilon _{x}Ax$ en un entier $n$ tel que $A(n)$ soit
vrai. Ackermann montre \`{a} la suite de la preuve de consistance proprement
dite qu'il est possible d'exhiber une fonction r\'{e}cursive qui donne une
borne pour les valeurs qui sont substitu\'{e}es aux $\epsilon $-termes des
formules. Ind\'{e}pendamment de son int\'{e}r\^{e}t intrins\`{e}que (en
particulier, elle fait appara\^{i}tre ''ce qui compte'' dans la preuve,
c'est-\`{a}-dire les arguments de cette fonction), l'existence et la nature
de cette fonction sont au coeur de la NCI de Kreisel. On expliquera donc en
d\'{e}tail la construction de cette fonction, en adaptant l'exposition \`{a}
l'utilisation de ce r\'{e}sultat par Kreisel, et on donnera ensuite les
\'{e}tapes de la preuve de ce que cette fonction borne effectivement les
valeurs des $\epsilon $-termes, en omettant la partie la plus fastidieuse et
la moins int\'{e}ressante de la d\'{e}monstration (S\"{a}tze 6 \`{a} 17 du
paragraphe 7 de \cite{acker}).

On peut poser le probl\`{e}me de la mani\`{e}re suivante : les valeurs
possibles pour les $\epsilon $-termes augmentent clairement avec le nombre
de substitutions que comporte la preuve. Une fois identifi\'{e}es les
variables pertinentes de la preuve, trouver une fonction de ces variables et
du nombre de ces substitutions qui borne les valeurs possibles, et une autre
qui borne \`{a} partir de ces seules variables le nombre de substitutions.
Si la premi\`{e}re partie de la t\^{a}che est facile \`{a} mener, la seconde
demandera davantage de travail. Avant cela, il nous faut d\'{e}finir un
codage des types d'ordres des ordinaux $<\epsilon _{0}$ sur les entiers qui
va permettre de manipuler au sein de l'arithm\'{e}tique habituelle les
indices de suite. Si cette op\'{e}ration peut s'expliquer par le souci
finitiste d'Ackermann, ce codage n'en reste pas moins facultatif, au sens
o\`{u} l'on pourrait apr\`{e}s tout manipuler directement les ordinaux. Mais
l'utilisation par Kreisel du r\'{e}sultat de bornage repose de mani\`{e}re
essentielle sur la possibilit\'{e} d'effectuer le calcul des bornes dans
l'arithm\'{e}tique-epsilon, de sorte qu'il est n\'{e}cessaire de
pr\'{e}senter pr\'{e}cis\'{e}ment les choses.

\subsubsection{Le codage des ordinaux}

Au vu de la preuve de terminaison, les ordinaux qu'on a besoin d'utiliser
sont tous des ordinaux de type $m$ pour un $m$ donn\'{e} au sens de la
d\'{e}finition suivante~:

\begin{definition}
Les ordinaux de la forme $\omega a+b$ sont de type $1$.

Si $\alpha _{1},\alpha _{2,}...\alpha _{i}$ sont des ordinaux de type $n$
tels que $\alpha _{1}>\alpha _{2}...>\alpha _{i},$ $\omega ^{\alpha
_{1}}+\omega ^{\alpha 2}...+\omega ^{\alpha _{i}}$ est un ordinal de type $%
n+1.$
\end{definition}

On cherche alors, pour chaque $m$, des relations d'ordre $<_{m}$ et une
bijection $f(\beta ,m)$\ qui associe \`{a} un ordinal $\beta $ de type $m$
un entier telles que $f$ soit un isomorphisme entre les ordinaux de type $m$
avec l'ordre habituel sur les ordinaux et les entiers \'{e}quip\'{e}s de $%
<_{m}.$

\begin{definition}
$f(\omega a+b,1)=$ $2^{a}(2b+1)-1$

$2^{a}(2b+1)-1<_{1}2^{c}(2d+1)-1$ si et seulement si $(a,b)<(c,d)\ $avec
l'ordre lexicographique

$f(\omega ^{\alpha _{1}}+\omega ^{\alpha _{2}}...+\omega ^{\alpha
_{i}},m+1)=2^{a_{1}}+2^{a_{2}}...+2^{a_{i}}$ o\`{u} $f(\alpha _{j},m)=a_{j}.$

$2^{a_{1}}+2^{a_{2}}...+2^{a_{i}}<_{m+1}2^{b_{1}}+2^{b_{2}}...+2^{b_{j}}$
o\`{u} $a_{1}>_{m}a_{2}>_{m}...>_{m}a_{i}$ et $%
b_{1}>_{m}b_{2}>_{m}...>_{m}b_{j}$ si et seulement si il existe un $k\leq
i,j $ tel que $a_{1}=b_{1},...a_{k-1}=b_{k-1}$ et $a_{k}<_{m}b_{k}$ ou si $%
j>i$ et $a_{k}=b_{k}$ pour tout $k\leq i.$
\end{definition}

Que $f$ et les $<_{m}$ r\'{e}alisent l'isomorphisme souhait\'{e} se
montrerait facilement par r\'{e}currence sur $m$\ en utilisant simplement
leur d\'{e}finition et l'unicit\'{e} \`{a} l'ordre fix\'{e} des
d\'{e}compositions propos\'{e}es des entiers.

On se donne pour la suite deux fonctions r\'{e}cursives primitives $\nu (x)$
et $\theta (x)$ telles que $\nu (2^{a}(2b+1)-1)=a$ et $\theta
(2^{a}(2b+1)-1)=b.$

On peut alors d\'{e}finir une classe de fonctions sp\'{e}cifiques, les
fonctions primitives r\'{e}cursives d'ordre fini, qui seront utilis\'{e}es
pour calculer la limite sup\'{e}rieure des valeurs des $\epsilon $-termes
clos. Ces fonctions sont construites de mani\`{e}re analogue aux fonctions
r\'{e}cursives primitives, moyennant une lib\'{e}ralisation du sch\'{e}ma de
r\'{e}currence. Soit $n$ un entier quelconque,\ si $g,h,\phi $ sont des
fonctions d\'{e}j\`{a} d\'{e}finies et que $\phi $ satisfait la
propri\'{e}t\'{e} $\phi (m)<_{n}m$ pour tout $m,$ on peut d\'{e}finir une
nouvelle fonction $f$ de la mani\`{e}re suivante :

$f(0,a)=g(0)$

$f(m,a)=h(a,m,f(\phi (m),a)$

On peut remarquer que le sch\'{e}ma de r\'{e}currence primitive d'ordre fini
ne fait d\'{e}pendre la valeur d'une fonction pour son argument que d'un
nombre fini de valeurs pr\'{e}c\'{e}demment calcul\'{e}es, ce qui le
distingue d'une d\'{e}finition par r\'{e}currence sur les ordinaux faisant
intervenir pour le calcul \`{a} la limite l'ensemble infini des valeurs
pr\'{e}c\'{e}dentes, et explique son acceptabilit\'{e} du point de vue
finitiste qui est celui d'Ackermann. Nous allons voir ensuite que les
fonctions r\'{e}cursives primitives d'ordre fini sont en fait des fonctions
r\'{e}cursives.

\subsubsection{Le bornage en fonction du nombre de substitutions}

Etant donn\'{e} l'ensemble de formules constitu\'{e}es par une preuve, on
d\'{e}signera \`{a} partir de maintenant par $m$ le degr\'{e} maximal d'un
terme dans la preuve, $e$ le nombre d'$\epsilon $-termes et $g$ le nombre de
cat\'{e}gories.

On d\'{e}finit par r\'{e}currence une fonction $\phi (m,a)$ qui donne la
valeur maximale des termes pour une substitution en fonction de la valeur
maximale $a$ d'un substituant pour les $\epsilon $-termes clos (il importe
ici de ne pas confondre valeur maximale d'un $\epsilon $-terme clos qui
d\'{e}pend directement des substituants aux cat\'{e}gories et valeur
maximale d'un terme qui d\'{e}pend ensuite des op\'{e}rations
\'{e}ventuellement appliqu\'{e}es aux $\epsilon $-termes clos) :

\begin{definition}
$\phi (0,a)=a$

$\phi (m+1,a)=[\phi (m,a)]^{2}+1$
\end{definition}

Il faut ici se souvenir de ce que les seuls symboles de fonction sont $%
^{\prime },d,+$ et $\times .$ L'id\'{e}e est qu'en un coup, on peut au plus
soit ajouter un, si la valeur maximale des $\epsilon $-termes clos est $\leq
1$, avec $^{\prime },$ soit multiplier le terme le plus grand par
lui-m\^{e}me sinon.

On d\'{e}finit ensuite par r\'{e}currence une fonction $\omega (m,n)$ qui
borne la valeur maximale des $\epsilon $-termes clos au bout de $n$
substitutions.

\begin{definition}
$\omega (m,0)=\phi (m,0)$

$\omega (m,n+1)=\phi (m,\omega (m,n))$
\end{definition}

Le cas de base d\'{e}coule de ce que la substitution initiale annule tous
les substituants. L'\'{e}tape de r\'{e}currence est justifi\'{e}e par le
fait que la valeur maximale d'un substituant pour $S_{n+1}$ est, telle que
celle-ci a \'{e}t\'{e} d\'{e}finie, inf\'{e}rieure ou \'{e}gale \`{a} la
valeur maximale d'un terme \`{a} l'\'{e}tape $S_{n}.$

\subsubsection{Le bornage du nombre de substitutions}

Reste donc \`{a} \'{e}valuer le nombre de substitutions. Cette
\'{e}valuation repose essentiellement sur la d\'{e}finition par
r\'{e}currence crois\'{e}e de deux fonctions : la premi\`{e}re $\tau
(c,p,n,a)$ calcule une borne sup\'{e}rieure pour l'indice d'une $p+1$%
-section qui commence par une substitution $S_{n}$ et une $p$-section
d'indice $a;$ la seconde $\kappa (c,p,n,a)$ calcule une borne sup\'{e}rieure
au sens de $<_{p}$\ pour une $p$-s\'{e}rie commen\c{c}ant en $S_{n}$ suivant
une $p$-s\'{e}rie d'indice $a.$

On suppose donn\'{e}e une fonction r\'{e}cursive $\eta (a,p)$ telle que $%
\eta (2^{a_{1}}+...+2^{a_{i}},p)=a_{1}$ o\`{u} les $a_{j}$ constituent une
suite strictement d\'{e}croissante.

On commence par donner deux d\'{e}finitions pr\'{e}liminaires, pour une
fonction $\psi (m,n,e)$ qui donne une borne pour le degr\'{e} d'une
substitution $S_{n}$ et pour une fonction $\lambda (a,p)$ qui calcule le
nombre de substitutions dont est compos\'{e}e une $p$-s\'{e}rie d'indice $a.$

\begin{definition}
$\psi (m,n,e)=2^{(\omega (m,n)+1)e}$
\end{definition}

En effet, le degr\'{e} est l'ordre relativement \`{a} un ensemble de
formules B(0,z$_{1}...z_{n}$), ...B(z,z$_{1}...z_{n}$). Or celles-ci
comporte au plus $(z+1)e$ termes, et $z\leq \omega (m,n+1)$ puisque $z$ est
la valeur d'un terme de $S_{n}.$

\begin{definition}
$\lambda (a,1)=1$

$\lambda (2^{a_{1}}+...+2^{a_{m}},p+1)=\lambda (a_{1},p)+...+\lambda
(a_{m},p)$
\end{definition}

La d\'{e}finition est cette fois imm\'{e}diatement claire.

On d\'{e}finit ensuite $\kappa (c,p,n,a)$ par cas

\begin{definition}
a) $\kappa (c,p,n,0)=0.$ Si $a\neq 0,$ on distingue

b) si $p=1,$ $\theta (a)\neq 0,$ $\kappa (c,p,n,a)=2^{\nu (a)}(2\times
\vartheta (a)-1)-1.$

c) si $p=1,$ $\theta (a)=0,$ $\kappa (c,p,n,a)=2^{\nu (a)-1}(2\times c+1)-1$

d) si $p>1$ et $a$ est pair, $\kappa (c,p,n,a)=a-1$

e) si $p>1$ et $a$ est impair, et $a=2^{a_{1}}-1,$

$\kappa (c,p,n,a)=\tau (c,p-1,n,\kappa (c,p-1,n,a_{1}))$

f) si $p>1$ et $a$ est impair, et $a=2^{a_{1}}+2^{a_{2}}+...+2^{a_{l}}$
o\`{u} les $a_{i}$ sont rang\'{e}s par ordre d\'{e}croissant

$\kappa (c,p,n,a)=2^{a_{1}}+\kappa (c,p,n+\lambda
(a_{1},p-1),2^{a_{2}}+...+2^{a_{l}}-1)$
\end{definition}

V\'{e}rifions informellement que $\kappa $ calcule bien ce qu'elle doit
calculer.\ Le cas a) est \'{e}vident. Les d\'{e}finitions pour b) et c)
majorent bien la 1-s\'{e}rie suivant une 1-s\'{e}rie d'indice $a.$ Pour le
cas b), ceci d\'{e}coule directement du fait que $\kappa (c,p,n,a)$ est le $%
<_{1}$ pr\'{e}d\'{e}cesseur imm\'{e}diat de $a$ par construction. Pour le
cas c), la valeur donn\'{e}e \`{a} c sera $\psi (m,n,e)$, en effet le
degr\'{e} de la substitution suivante sera $\leq \psi (m,n,e)$ comme la
valeur maximale possible pour un terme n'a pas augment\'{e} d'une
substitution \`{a} l'autre. d) sert juste \`{a} assurer que $\kappa
(c,p,n,a)<_{p}a$ pour le bon fonctionnement de la r\'{e}currence; en
pratique cela n'intervient pas, car la seule substitution d'indice 0 est la
substitution initiale, et que par cons\'{e}quent la seule $p$-s\'{e}rie
d'indice 0 est la 1-s\'{e}rie qu'elle constitue, de sorte que les indices
des autres sections sont toujours impaires. e) correspond au cas o\`{u} la $p
$-s\'{e}rie P est compos\'{e}e d'une unique $p-1$-s\'{e}rie M. Par
cons\'{e}quent, la $p$-s\'{e}rie suivant P commence par une $p-1$-s\'{e}rie
qui succ\`{e}de imm\'{e}diatement \`{a} M, de sorte qu'on peut l'\'{e}valuer
\`{a} l'aide de la fonction $\tau $ \`{a} laquelle on fournit comme majorant
de l'indice de M $\kappa (c,p-1,n,a).$ f) correspond \`{a} la situation
o\`{u} P est compos\'{e}e de $l$ $p-1$-s\'{e}ries, et fait avancer le calcul
en majorant au moyen de l'indice de la premi\`{e}re de ces s\'{e}ries et de $%
\kappa $ appliqu\'{e}e \`{a} la somme des indices des $p-1$-s\'{e}ries
restantes, de sorte que la valeur de l'argument d\'{e}cro\^{i}t comme il
convient, tout en tenant compte du d\'{e}calage de $S_{n}$ \`{a} $%
S_{n+\lambda (a_{1},p-1)}.$

\begin{definition}
$\tau (c,p,n,0)=0$

$\tau (c,p,n,a)=2^{a}+\tau (c,p,n+\lambda (a,p),\kappa (c,p,n+\lambda
(a,p),a))$
\end{definition}

Le cas de base est \'{e}vident. Pour l'\'{e}tape de r\'{e}currence, $\tau $
doit calculer l'indice $i$ d'une $p+1$-s\'{e}rie P commen\c{c}ant par une $p$%
-s\'{e}rie M d'indice $a.$ On a $i=2^{a}+i^{\prime }$ ou $i^{\prime }$ doit
\^{e}tre l'indice correspondant au reste des $p$-s\'{e}ries qui composent P,
vues comme une $p+1$-s\'{e}rie. C'est pr\'{e}cis\'{e}ment ce que calcul $%
\tau $ si on l'applique \`{a} un majorant de l'indice de la $p$-s\'{e}rie
qui suit M, ce qui nous est pr\'{e}cis\'{e}ment donn\'{e} par $\kappa
(c,p,n+\lambda (a,p),a),$ en tenant compte du d\'{e}calage de $S_{n}$ \`{a} $%
S_{n+\lambda (a_{1},p)}$. On remarque que la d\'{e}finition de $\tau $ fait
intervenir essentiellement une r\'{e}currence d'ordre sup\'{e}rieur; on n'a
en effet aucune raison d'avoir $\kappa (c,p,n+\lambda (a,p),a)<_{0}a.$ Par
contre, pour la bonne d\'{e}finition de $\tau ,$ il faut maintenant
v\'{e}rifier que $\kappa (c,p,n,a)<_{p}a$ quand $a\neq 0.$

La d\'{e}monstration se fait par r\'{e}currence sur $p.$ On a vu que les cas
a), b) et d) ne posaient pas de probl\`{e}me. Pour le cas c), il faut en
fait voir la variable c comme un param\`{e}tre : $\kappa $ est bien
d\'{e}finie pour des valeurs convenables de c, en particulier lorsque $%
c=\psi (m,n,e).$ Supposons que $\kappa (c,p,n,a)<_{p}a$ quand $a\neq 0$ et
montrons la propri\'{e}t\'{e} pour $p+1.$ Par l'hypoth\`{e}se de
r\'{e}currence, $\tau $ est calculable pour des $p^{\prime }\leq p$ et donc $%
\kappa $ pour des $p^{\prime }\leq p+1.$ On d\'{e}montre ensuite un lemme
par r\'{e}currence transfinie imbriqu\'{e}e dans la pr\'{e}c\'{e}dente que $%
\eta (\tau (c,p,n,a),p)=a.$ Le cas $a=0$ est imm\'{e}diat. Supposons
maintenant que la propri\'{e}t\'{e} tient pour des $a^{\prime }<_{p}a.$

(1) $\kappa (c,p,n+\lambda (a,p),a)<_{p}a$ par hypoth\`{e}se de
r\'{e}currence

(2) $\eta (\tau (c,p,n+\lambda (a,p),\kappa (c,p,n+\lambda
(a,p),a),p)=\kappa (c,p,n+\lambda (a,p),a)$ puisque (1) permet d'appliquer
l'hypoth\`{e}se de r\'{e}currence transfinie.

(3) $\eta (\tau (c,p,n+\lambda (a,p),\kappa (c,p,n+\lambda (a,p),a),p)<_{p}a$
par (1) +(2)

(4) si $\eta (b,p)<_{p}a,$ alors $\eta (2^{a}+b,p)=a$ par d\'{e}finition de $%
\eta .$

(5) $\eta (2^{a}+\tau (c,p,n+\lambda (a,p),\kappa (c,p,n+\lambda
(a,p),a),p)=a$ en appliquant (4) \`{a} (3) avec $b=\tau (c,p,n+\lambda
(a,p),\kappa (c,p,n+\lambda (a,p),a)$

(6) $\eta (\tau (c,p+1,n,a),p+1)=a$ par d\'{e}finition de $\tau .$

Examinons maintenant le cas e). On veut $\kappa (c,p+1,n,a)<_{p+1}a$ o\`{u} $%
a$ est $2^{a_{1}}-1.$

(1) $\kappa (c,p+1,n,a)=\tau (c,n,\kappa (c,p,n,a_{1}))$ par d\'{e}finition
de $\kappa $

(2) $\eta (\tau (c,p,n,\kappa (c,p,n,a_{1}),p)=\kappa (c,p,n,a_{1})$ par le
lemme

(3) $\kappa (c,p,n,a_{1})<_{p}a_{1}$ par l'hypoth\`{e}se de r\'{e}currence

(4) $\eta (a,p)=a_{1}$ par hypoth\`{e}se sur $a.$

(5) $\eta (\kappa (c,p+1,n,a),p)<_{p}\eta (a,p)$ par (1), (2), (3) et (4)

(6) si $\eta (a,p)<_{p}\eta (b,p),$ alors $a<_{p+1}b$ par d\'{e}finition de $%
\eta .$

(7) $\kappa (c,p+1,n,a)<_{p+1}a$ en appliquant (6) \`{a} (5).

Reste le cas f).

? $\kappa (c,p+1,n,a)<_{p+1}2^{a_{1}}+2^{a_{2}}+...+2^{a_{l}}-1$

? $2^{a_{1}}+\kappa (c,p+1,n+\lambda
(a_{1},p),2^{a_{2}}+...+2^{a_{l}}-1)<2^{a_{1}}+2^{a_{2}}+...+2^{a_{l}}-1$
par d\'{e}finition de $\kappa $

? $\kappa (c,p+1,n+m,2^{a_{l}}-1)<2^{a_{l}}$ par application
r\'{e}p\'{e}t\'{e}e de la d\'{e}finition de $\kappa ,$ pour un certain $m$,
ce qui nous ram\`{e}ne au cas e).

On se contente pour finir du sch\'{e}ma de la d\'{e}monstration en omettant
de d\'{e}montrer deux propositions. Les d\'{e}monstrations omises se font
presque toutes sur le principe d'une induction sur $p$ o\`{u} interviennent
simultan\'{e}ment $\tau $ et $\kappa ,$ mais plusieurs \'{e}tapes sont
n\'{e}cessaires avant de pouvoir d\'{e}montrer les propri\'{e}t\'{e}s de $%
\lambda $ et $\tau $ dont on a besoin. On d\'{e}finit un pr\'{e}dicat T
\`{a} quatre arguments qui s'appliquera \`{a} $(\psi (m,n,e),p,n,a)$ si $a$
est l'indice d'une $p$-s\'{e}rie commen\c{c}ant par $S_{n}$ pour une preuve
dont les constantes sont $m$ et $e$. On montre alors premi\`{e}rement que le
majorant donn\'{e} par $\tau $ se comporte bien avec $\lambda $ et
deuxi\`{e}mement que $\tau $ et $\lambda $ ensemble se comportent bien
vis-\`{a}-vis de $\leq _{p}.$

\begin{proposition}
Si $T(\psi (m,n,e),p+1,n,a)$, alors

$\lambda (a,p+1)\leq \lambda (\tau (\psi (m,n,e),p,n,\eta (a,p)),p+1)$
\end{proposition}

\begin{proposition}
Si $a\leq _{p}b$ et $T(\psi (m,n,e),n,a),$ alors

$\lambda (\tau (\psi (m,n,e),p,n,a),p+1)\leq \lambda (\tau (\psi
(m,n,e),p,n,b),p+1)$
\end{proposition}

Supposons donc qu'on connaisse les indices $a_{1}...a_{t}$ des $t$ $g$%
-s\'{e}ries qui composent l'unique $g+1$-s\'{e}ries que constitue la suite
compl\`{e}te des substitutions, la proposition [ ] nous dit que

$\lambda (2^{a_{1}}+...+2^{a_{t}},g+1)\leq \lambda (\tau (\psi
(m,n,e),g,1,a_{1}),g+1)$

L'\'{e}valuation du nombre de substitutions ne d\'{e}pend alors plus de $t$%
.La seconde proposition nous dit qu'il suffit de trouver un majorant $b$\ au
sens de $\leq _{g}$ de $a_{1}$ pour avoir

$\lambda (2^{a_{1}}+...+2^{a_{t}},g+1)\leq \lambda (\tau (\psi
(m,n,e),g,1,b),g+1)$

On d\'{e}finit enfin une fonction $\rho (n,e)$.

\begin{definition}
$\rho (1,e)=2^{e+1}-1$

$\rho (n+1,e)=2^{\rho (n,e)}-1$
\end{definition}

On voit que $\rho (1,e)$ majore au sens de $<_{1}$ l'indice d'une $1$%
-s\'{e}rie initiale, puisque $e$ est le nombre d'$\epsilon $-termes clos de
la preuve et que si $\rho (n,e)$ majore au sens de $<_{n}$ l'indice d'une $n$%
-s\'{e}rie initiale, alors $\rho (n+1,e)$ majore au sens de $<_{n+1}$
l'indice d'une $n+1$-s\'{e}rie initiale. Donc $a_{1}<_{g}\rho (g,e).$ On en
d\'{e}duit

$\lambda (2^{a_{1}}+...+2^{a_{t}},g+1)\leq \lambda (\tau (\psi
(m,n,e),g,1,\rho (g,e)),g+1)$ o\`{u} la seconde expression ne d\'{e}pend
plus des indices des s\'{e}ries de la suite des substitutions tels qu'on
pourrait les calculer en appliquant effectivement la m\'{e}thode de
substitution. D'o\`{u} le th\'{e}or\`{e}me final.

\begin{theorem}
Soit une preuve de l'arithm\'{e}tique-epsilon comportant $e$ $\epsilon $%
-termes correspondant \`{a} $g$ cat\'{e}gories distinctes et telle que le
degr\'{e} des termes ne d\'{e}passe pas $m,$ la valeur maximale des $%
\epsilon $-termes clos dans la substitution finale est $born(m,e,g)=\omega
(m,\lambda (\tau (\psi (m,n,e),g,1,\rho (g,e)),g+1)).$
\end{theorem}

\section{La No-counterexample interpretation}

\subsection{La notion d'interpr\'{e}tation}

Kreisel d\'{e}finit une notion g\'{e}n\'{e}rale d'interpr\'{e}tation qui
convient \`{a} un syst\`{e}me d\'{e}ductif quelconque pour
l'arithm\'{e}tique. La sp\'{e}cificit\'{e} de la d\'{e}marche est d'exiger
que l'interpr\'{e}tation du syst\`{e}me se fasse \`{a} l'int\'{e}rieur de
celui-ci afin de contr\^{o}ler que l'interpr\'{e}tation elle-m\^{e}me est
bien finitiste. On se donne un codage r\'{e}cursif des formules du
syst\`{e}me.

\begin{definition}
Une formule $A[x_{1}...x_{n},f_{1}...f_{m}]$ sans variables li\'{e}es, et
dont les variables libres sont parmi les variables d'individus $%
x_{1}...x_{n} $ et les variables de fonctions $f_{1}...f_{m}$ est
v\'{e}rifiable si, pour tous num\'{e}raux $a_{1}...a_{n}$ et toutes
fonctions r\'{e}cursives $\phi _{1}...\phi _{m},$ $%
A[x_{i}:=a_{i},f_{j}:=g_{j}]$ est vraie.
\end{definition}

\begin{definition}
On appelle alors interpr\'{e}tation d'un syst\`{e}me $\Sigma $ une fonction
calculable $f(n,a)$ telle que :

$\alpha )$ $f(n,a)$ est le num\'{e}ro d'une formule $A_{n}$ dont toutes les
variables sont libres lorsque $a$ est le num\'{e}ro d'une formule $A$ de $%
\Sigma .$

$\beta )$ Si on a une preuve de $A$ dans $\Sigma ,$ on peut trouver \`{a}
partir de la preuve un $n$ tel que $A_{n}$ est v\'{e}rifiable.

$\gamma )$ Si on a une preuve de $\sim A$ dans $\Sigma ,$ pour chaque $n,$
on trouve une instanciation des variables libres de $A_{n}$ qui rend $A_{n}$
fausse.

$\delta )$ Si on a une preuve de $A\rightarrow B$ dans $\Sigma ,$ on trouve
une fonction calculable $g(n)$ telle que si $A_{n}$ est v\'{e}rifiable $%
B_{g(n)}$ l'est aussi.
\end{definition}

Le but est bien de capturer le ''sens finitiste'' des preuves
non-constructives, c'est-\`{a}-dire de passer de preuves d'\'{e}nonc\'{e}s
quantifi\'{e}s \`{a} des formules ayant un contenu num\'{e}rique
imm\'{e}diat. Par exemple, si $\forall xAx$ est une formule universelle, $Ax$
sera une bonne mani\`{e}re d'interpr\'{e}ter une preuve quelconque de $%
\forall xAx,$ car \`{a} partir de cette preuve, on peut tirer une preuve de $%
An$ pour n'importe quel num\'{e}ral $n.$ Mais on a vu dans l'introduction
que l'interpr\'{e}tation ne saurait \^{e}tre aussi directe. Inversement, il
est toujours possible de trivialiser la notion d'interpr\'{e}tation. Par
exemple, si $e(a)$ est le code de $\sim A$ quand $a$ est le code de $A$ et
si $Prov_{\Sigma }(x,y)$ est un pr\'{e}dicat de prouvabilit\'{e} pour $%
\Sigma ,$ on peut montrer facilement que, sous hypoth\`{e}se de la
consistance de $\Sigma $, l'association \`{a} $A\ $de $\sim Prov_{\Sigma
}(x,e(a))$ constitue une interpr\'{e}tation pour $\Sigma .$ On veut au
contraire que l'interpr\'{e}tation d'une formule conserve autant que
possible la \textit{signification }de celle-ci.

La NCI est donn\'{e}e pour une extension de PA$_{\epsilon }$ dans laquelle
on ajoute des symboles de variables de fonctions, que l'on ne quantifie pas.
Pr\'{e}sentons d'abord le principe de la NCI, avant de voir en quoi elle
peut constituer une interpr\'{e}tation de PA$_{\epsilon }.$ Soit $A$ une
formule en forme pr\'{e}nexe

$\forall x_{1}\exists y_{1}...\forall x_{n}\exists y_{n}A^{\prime
}(x_{1}...x_{n},y_{1}...y_{n})$

o\`{u} $A^{\prime }$ est sans quantificateur. L'interpr\'{e}tation
na\"{i}ve, dont on a vu qu'elle ne marchait pas, consistait \`{a} chercher
des fonctions calculables $\phi _{1}...\phi _{n}$ telles que

$A^{\prime }(x_{1}...x_{n},\phi _{1}(x_{1})...\phi _{n}(x_{1}...x_{n}))$

soit v\'{e}rifiable. Dans les termes de Kreisel, ceci constituerait une
\textit{Erf\"{u}llung }de $A.$ Consid\'{e}rons ce que serait une \textit{%
Erf\"{u}llung }de $\sim A.$ $\sim A$ elle-m\^{e}me est \'{e}quivalente \`{a}
:

$\exists x_{1}\forall y_{1}...\exists x_{n}\forall y_{n}\sim A^{\prime
}(x_{1}...x_{n},y_{1}...y_{n})$

Donc l'\textit{Erf\"{u}llung }serait donn\'{e}e par des fonctions
calculables $\psi _{1}...\psi _{n}$ (o\`{u} $\psi _{1}$ est une fonction
0-aire) telles que

$\sim A^{\prime }(\psi _{1}...\psi _{n}(y_{1}...y_{n-1}),y_{1}...y_{n})$

soit v\'{e}rifiable. Ceci constituerait un contre-exemple \`{a} $A$ au sens
o\`{u} les fonctions donnent pour toute instanciation $b_{1}...b_{n}$ des $%
y_{1}...y_{n}\ $des valeurs $\psi _{1}...\psi _{n}(b_{1}...b_{n-1})$ qui
montrent que les $b_{1}...b_{n}$ ne sont pas de bons t\'{e}moins pour les
existentiels de la formule de d\'{e}part. Faute de pouvoir obtenir les
fonctions calculables $\phi _{1}...\phi _{n},$ l'id\'{e}e est de montrer
qu'il n'y a pas d'\textit{Erf\"{u}llung }de $\sim A,$ autrement dit pas de
contre-exemple \`{a} $A.$ Pour chaque $\psi _{1}...\psi _{n},$ on veut donc
trouver des $b_{1}...b_{n}$ tels que $A^{\prime }(\psi _{1}...\psi
_{n}(b_{1}...b_{n-1}),b_{1}...b_{n}).$ Ceci revient \`{a} exhiber des
fonctionnels $\chi _{1}...\chi _{n}$ dont les variables libres sont parmi
les variables de fonctions $f_{1}...f_{n}$ telles que

$A^{\prime }(f_{1}...f_{n}(\chi _{1}...\chi _{n-1}),\chi _{1}...\chi _{n})$

soit v\'{e}rifiable.

Si on arrive \`{a} passer d'une preuve de $A$ \`{a} l'existence d'une suite
de tels fonctionnels, ceci nous donnera bien une interpr\'{e}tation,
moyennant une \'{e}num\'{e}ration ad\'{e}quate des fonctionnels en question (%
$A_{m}$ sera alors la formule $A^{\prime }(f_{1}...f_{n}(\chi _{1}...\chi
_{n-1}),\chi _{1}...\chi _{n})$ o\`{u} $\chi _{1}...\chi _{n}$ est la $m$%
-i\`{e}me suite de fonctionnels) et, \textit{last but not least}, la
v\'{e}rification des conditions $\gamma )$ et $\delta )$.
L'interpr\'{e}tation se fera donc dans un formalisme comportant des
variables libres de fonctions.

\subsection{V\'{e}rification de la condition $\protect\beta )$}

Voyons comment la preuve d'Ackermann nous permet de v\'{e}rifier
relativement facilement la condition $\beta )$. On commence par donner
l'id\'{e}e g\'{e}n\'{e}rale qui est de se servir de symboles de variables de
fonctions pour param\'{e}trer la preuve d'Ackermann et le bornage qui en
r\'{e}sulte. On part de

$\vdash _{PA_{\epsilon }}\forall x_{1}\exists y_{1}...\forall x_{n}\exists
y_{n}A^{\prime }(x_{1}...x_{n},y_{1}...y_{n})$

Il en d\'{e}coule

$\vdash _{PA_{\epsilon }^{\prime }}\exists y_{1}...\exists y_{n}A^{\prime
}(f_{1}...f_{n}(y_{1}...y_{n-1}),y_{1}...y_{n})$

o\`{u} $PA_{\epsilon }^{\prime }$ ne diff\`{e}re de $PA_{\epsilon }$ que par
l'adjonction aux langages de symboles de variables libres de fonctions comme
$f_{1}...f_{n}$. En d\'{e}finissant une suite ad\'{e}quate d'$\epsilon $%
-termes

$e_{y_{n}}=\epsilon _{y_{n}}A^{\prime
}(f_{1}...f_{n}(y_{1}...y_{n-1}),y_{1}...y_{n})$

$e_{y_{k-1}}=\epsilon _{y_{k-1}}A^{\prime
}(f_{1}...f_{n}(y_{1}...y_{k-1}e_{k...}e_{n-1}),y_{1}...y_{k-1}e_{k...}e_{n})
$

on obtient une formule \'{e}quivalente \`{a} la pr\'{e}c\'{e}dente o\`{u}
les $\epsilon $-termes ont remplac\'{e} les quantificateurs existentiels.

$\vdash _{PA_{\epsilon }^{\prime }}A^{\prime
}(f_{1}...f_{n}(e_{1}...e_{n-1}),e_{1}...e_{n})$

Supposons un instant qu'on ait choisi d'instancier les $x_{i}$ par des
fonctions r\'{e}cursives d\'{e}termin\'{e}es $\psi _{1}...\psi _{n}$, de
sorte qu'on ait

$\vdash _{PA_{\epsilon }^{\prime \prime }}A^{\prime }(\psi _{1}...\psi
_{n}(e_{1}^{\prime }...e_{n-1}^{\prime }),e_{1}^{\prime }...e_{n}^{\prime })$

o\`{u} $e_{j}^{\prime }=e_{j}[\overrightarrow{f}:=\overrightarrow{\psi }]\ $%
et $PA_{\epsilon }^{\prime \prime }$ ne diff\`{e}re de $PA_{\epsilon }$ que
par l'adjonction des symboles de fonctions $\psi _{1}...\psi _{n}$. La
m\'{e}thode de substitution va s'appliquer \`{a} $PA_{\epsilon }^{\prime
\prime }$ comme \`{a} $PA_{\epsilon },$ de sorte qu'\`{a} partir de la
preuve pr\'{e}c\'{e}dente, on obtient des num\'{e}raux $m_{1}...m_{n}$ pour
lesquels

$A^{\prime }(\psi _{1}...\psi _{n}(m_{1}....m_{n-1}),m_{1}...m_{n})$ (*)

soit vrai. De plus, on peut donner \`{a} l'avance une limite sup\'{e}rieure $%
m$ pour les $m_{1}...m_{n}$ en utilisant une fonction analogue \`{a} $%
born(m,e,g)$. $born(m,e,g)$ ne convient pas directement, car on a enrichi le
langage de $PA_{\epsilon }$ avec de nouveaux symboles de fonctions, de sorte
que la fonction $\phi (m,a)$ doit \^{e}tre remplac\'{e}e par une fonction $%
\phi ^{\prime }(m,a)$ qui majore les valeurs des termes en tenant compte des
nouvelles fonctions.

Il appara\^{i}t ainsi que $born(m,e,g)$ d\'{e}pend des fonctions qu'on
substitue aux symboles de variables de fonctions. L'id\'{e}e de Kreisel est
alors de capturer pr\'{e}cis\'{e}ment cette d\'{e}pendance en rempla\c{c}ant
$born(m,e,g)$ par une fonctionnelle $born^{\prime }[f_{1}...f_{n}]$ (les
param\`{e}tres $m,e,g$ \'{e}tant fournis par la preuve)\ comportant comme
variables libres de fonctions les $f_{1}...f_{n}.$ Pour coller exactement
avec la notion d'interpr\'{e}tation qui a \'{e}t\'{e} d\'{e}finie, on pourra
alors exprimer l'interpr\'{e}tation d'une formule au moyen d'un sch\'{e}ma
de minimisation.

La d\'{e}monstration de $\beta )$ est pour ainsi dire faite. Mais pour faire
les choses convenablement, il faut d\'{e}finir la classe de fonctionnels.
Avant cela, on commence par d\'{e}finir des avatars $\phi _{\overrightarrow{x%
}}^{\prime }$ et $\omega _{m,\overrightarrow{x}}^{\prime }$ des fonctions $%
\phi \ $et $\omega $ avec les m\^{e}mes significations intuitives$.$ Le
param\`{e}tre$\ m$ correspond toujours au degr\'{e} maximal des termes. On
d\'{e}signe par $g_{j},$ $j\leq l,$ les $l$ symboles de fonctions -
constantes ou variables - $k_{j}$-aires qui apparaissent dans la preuve.

\begin{definition}
$\phi _{\overrightarrow{x}}^{\prime }[g_{j},a,1]=max_{j\leq
l,a_{1}...a_{k_{j}}\leq a}(g_{j}(a_{1}...a_{k_{j}}))$

$\phi _{\overrightarrow{x}}^{\prime }[g_{j},a,n+1]=\phi _{\overrightarrow{x}%
}^{\prime }[g_{j},\phi _{\overrightarrow{x}}^{\prime }[g_{j},a,n],1],$

$\omega _{m,\overrightarrow{x}}^{\prime }(0)=\phi _{\overrightarrow{x}%
}^{\prime }[g_{j},0,m]$

$\omega _{m,\overrightarrow{x}}^{\prime }(n+1)=\phi _{\overrightarrow{x}%
}^{\prime }[g_{j},\omega _{m,\overrightarrow{x}}^{\prime }(n),m]$
\end{definition}

$max_{j\leq l,a_{1}...a_{k_{j}}\leq a}(g_{j}(a_{1}...a_{k_{j}}))$ est
repr\'{e}sentable dans PA$_{\epsilon }$ par $\epsilon _{x}(\forall
\overrightarrow{y}\leq a,g_{1}(\overrightarrow{y})\leq x\wedge ...\wedge
g_{l}(\overrightarrow{y})\leq x),$ donc $\phi _{\overrightarrow{x}}^{\prime
} $ d\'{e}finie par r\'{e}currence \`{a} partir de cette fonction
\'{e}galement, et $\omega _{m,\overrightarrow{x}}^{\prime }$ du m\^{e}me
coup.

On modifie la d\'{e}finition de $\psi $ en $\psi ^{\prime }$en rempla\c{c}%
ant $\omega $ par $\omega _{m,\overrightarrow{x}}^{\prime }$

\begin{definition}
$\psi _{m}^{\prime }(n,e)=2^{(\omega _{m,\overrightarrow{x}}^{\prime
}(n)+1)e}$
\end{definition}

\begin{definition}
L'ensemble des fonctionnelles primitives r\'{e}cursives d'ordre fini
contenant les variables d'individus $a_{i}$ et les variables de fonctions $%
f_{j}$ est le plus petit ensemble E tel que

a) les fonctions primitives r\'{e}cursives d'ordre fini de variables $a_{i}$
sont dans E

b) si $g(x,y)$ est une fonction primitive r\'{e}cursive d'ordre fini et $%
\chi \in E,$ $g(a_{i},\chi )\in E.$

c) si $\chi \in E,\ \omega _{m,\overrightarrow{x}}^{\prime }(\chi )$ est
dans E.

d) si $\chi _{1}...\chi _{n}\in E$ et $f_{j}$ est n-aire, $f_{j}$
appliqu\'{e}e n-\'{e}l\'{e}ments parmi \{$\chi _{1}...\chi
_{n},a_{1}...a_{n} $\}est dans E.

e) si $\chi \in E$ et $A(x,...a_{i}...,...f_{j}...)$ est une formule sans
quantificateur dont les variables libres sont $x$ et des variables parmi les
$a_{i},f_{j}$ alors $\epsilon _{x}x\leq \chi \wedge
A(x,...a_{i}...,...f_{j}...)\in E.$
\end{definition}

La d\'{e}finition de Kreisel \cite{kreis2} est corrig\'{e}e conform\'{e}ment
\`{a} Kreisel \cite{kreis3} : l'ajout d'une clause assurant la cl\^{o}ture
par minimisation born\'{e}e est n\'{e}cessaire pour rejoindre la
d\'{e}finition pr\'{e}cise d'une interpr\'{e}tation. Afin d'\'{e}viter dans
un premier temps les codages qui compliqueraient l'\'{e}tablissement du
r\'{e}sultat (en obligeant \`{a} tenir compte des fonctions de codage pour
l'accroissement des valeurs des termes), on a en outre modifi\'{e} la
d\'{e}finition de mani\`{e}re \`{a} y faire entrer des fonctionnelles \`{a}
plusieurs variables de fonctions d'arit\'{e}s quelconques. On peut
maintenant \'{e}noncer le th\'{e}or\`{e}me.

\begin{theorem}
\label{100}Si on a une preuve de $\forall x_{1}\exists y_{1}...\forall
x_{n}\exists y_{n}A^{\prime }(x_{1}...x_{n},y_{1}...y_{n})$ dans PA$%
_{\epsilon },$ alors il existe des fonctionnelles r\'{e}cursives d'ordre
fini $\eta _{1}...\eta _{n}$ dont les variables libres sont parmi $%
f_{1}...f_{n}$ telles que $A^{\prime }(f_{1}...f_{n}(\eta _{1}...\eta
_{n-1}),\eta _{1}...\eta _{n})$ est v\'{e}rifiable.
\end{theorem}

D'une part, on observe que les fonctionnelles qui correspondent aux
fonctions du th\'{e}or\`{e}me de bornage sont de l'esp\`{e}ce voulue $\omega
_{m,\overrightarrow{x}}^{\prime }(n)$ est une fonctionnelle r\'{e}cursive
primitive, donc $\psi _{m}^{\prime }(n,e)$ aussi. $\lambda $ et $\rho $ sont
des fonctions primitives r\'{e}cursives, $\tau $ et $\kappa $ sont des
fonctions primitives r\'{e}cursives d'ordre fini, donc $\lambda (\tau (\psi
_{m}^{\prime }(n,e),g,1,\rho (g,e)),g+1)$ est une fonctionnelle primitive
r\'{e}cursive d'ordre fini, donc $\omega _{m,\overrightarrow{x}}^{\prime
}(\lambda (\tau (\psi _{m}^{\prime }(n,e),g,1,\rho (g,e)),g+1))$
\'{e}galement. On note $born_{m,n,e}^{\prime }[f_{1}...f_{n}]$cette
derni\`{e}re fonctionnelle param\'{e}tr\'{e}e par $m,n,e.$

D'autre part, on constate que le th\'{e}or\`{e}me de bornage reste valable
(il suffit pour cela de v\'{e}rifier qu'avec les changements effectu\'{e}s
sur $\omega ,\phi $ et $\psi $ la d\'{e}monstration reste correcte). Etant
donn\'{e}e une preuve de l'arithm\'{e}tique-epsilon dans le langage
augment\'{e}e des fonctions r\'{e}cursives $\overrightarrow{\psi },$
comportant $e$ $\epsilon $-termes correspondant \`{a} $g$ cat\'{e}gories
distinctes et telle que le degr\'{e} des termes ne d\'{e}passe pas $m,$ la
valeur maximale des $\epsilon $-termes clos dans la substitution finale est $%
born_{m,n,e}^{\prime }[f_{1}:=\psi _{1}...f_{n}:=\psi _{n}].$

Soient $\psi _{1}...\psi _{n}$ des fonctions r\'{e}cursives susceptibles
d'instancier les variables de fonctions $f_{1}...f_{n},$ on obtient \`{a}
partir d'une preuve $\pi $ de $A$ une preuve $\pi ^{\prime }(m,n,e)$ de $%
A^{\prime }(\psi _{1}...\psi _{n}(e_{1}^{\prime }...e_{n-1}^{\prime
}),e_{1}^{\prime }...e_{n}^{\prime }).$ On peut choisir les $\pi ^{\prime }$
de mani\`{e}re \`{a} ce que les triplets $(m,n,e)$ soient les m\^{e}mes,
ind\'{e}pendamment du choix des $\psi _{1}...\psi _{n}$ (ils d\'{e}pendent
seulement des constantes de la preuve $\pi $).\ On se donne un codage $%
\alpha _{n}$ des n-uplets et des fonctions $p_{n}^{1}...p_{n}^{n}$ de
d\'{e}codage, tout cela r\'{e}cursif primitif. On note $\beta $
fonctionnelle $\alpha _{n}(born_{m,n,e}^{\prime },...,born_{m,n,e}^{\prime
}).$ On d\'{e}finit alors $\eta _{i}=p_{n}^{i}(\epsilon _{z}(z\leq \beta
\wedge A^{\prime
}(f_{1}...f_{n}(p_{n}^{1}(z)...p_{n}^{n}(z)),p_{n}^{1}(z)...p_{n}^{n}(z))).$
Les $\eta _{i}$ sont bien \`{a} leur tour des fonctionnelles r\'{e}cursives
primitives d'ordre fini $.$ Le th\'{e}or\`{e}me 15 et le th\'{e}or\`{e}me de
bornage adapt\'{e} nous garantissent l'existence de $m_{1}...m_{n}\leq
born_{m,n,e}^{\prime }[f_{1}:=\psi _{1}...f_{n}:=\psi _{n}]$ pour lesquels $%
A^{\prime }(\psi _{1}...\psi _{n}(m_{1}....m_{n-1}),m_{1}...m_{n})$ est
vraie.

Par cons\'{e}quent, $A^{\prime }(f_{1}...f_{n}(\eta _{1}...\eta _{n-1}),\eta
_{1}...\eta _{n})$ est v\'{e}rifiable.

\subsection{Modularit\'{e} de l'interpr\'{e}tation}

Les conditions $\gamma )$ et $\delta )$ reviennent \`{a} faire peser des
conditions de modularit\'{e} sur l'interpr\'{e}tation propos\'{e}e. On veut
pouvoir induire de ce qu'il existe une interpr\'{e}tation pour $\sim $A
qu'il n'existe pas d'interpr\'{e}tation pour A; et on veut passer d'une
interpr\'{e}tation pour A \`{a} une interpr\'{e}tation pour B lorsqu'on
poss\`{e}de une preuve de $A\rightarrow B.$ Ces conditions
suppl\'{e}mentaires sont non triviales (voir \cite{kohl}); elles distinguent
la NCI de Kreisel des versions que l'on obtient via traduction de G\"{o}del
et interpr\'{e}tation fonctionnelle (c'est-\`{a}-dire les versions
tir\'{e}es de \cite{godel}).

Voyons d'abord o\`{u} se loge la difficult\'{e} dans la preuve de $\gamma )$%
. L'id\'{e}e de la preuve est simple : r\'{e}appliquer la m\'{e}thode de
substitution \`{a} une variante de $\sim A^{\prime }$ contenant les
fonctionnelles \`{a} r\'{e}futer ce qui nous fournira parmi les substituants
finaux les fonctions r\'{e}cursives recherch\'{e}es. Mais la r\'{e}alisation
de l'id\'{e}e est un peu d\'{e}licate. Supposons qu'on ait une preuve dans PA%
$_{\epsilon }$ de $\sim A$ o\`{u} $A$ est $\forall x_{1}\exists
y_{1}...\forall x_{n}\exists y_{n}A^{\prime }(x_{1}...x_{n},y_{1}...y_{n})$
et que les fonctionnelles \`{a} r\'{e}futer $\eta _{1}...\eta _{n}$ dont les
variables soient repr\'{e}sentables dans une extension de PA$_{\epsilon }$
par des $\epsilon $-termes $t_{1}...t_{n}$ (la notion de repr\'{e}sentation
est la m\^{e}me que pour les fonctions, cela signifie que l'on peut prouver $%
\eta _{i}=t_{i}$).

On d\'{e}finit une suite d'$\epsilon $-termes $e_{1}...e_{n}$ par :

$e_{1}=\epsilon _{x_{1}}\forall y_{1}...\exists x_{n}\forall y_{n}\sim
A^{\prime }(x_{1}...x_{n},y_{1}...y_{n})$

$e_{k+1}=\epsilon _{x_{k+1}}\forall y_{k+1}...\exists x_{n}\forall y_{n}\sim
A^{\prime
}(e_{1}...e_{k}(a_{1}...a_{k-1}),x_{k+1}...x_{n},a_{1}...a_{k},y_{k+1}...y_{n})
$

Il est clair qu'\`{a} partir d'une preuve de $\sim A$ on peut obtenir une
preuve de $\sim A^{\prime }(e_{1}...e_{n}(a_{1}...a_{n-1}),a_{1}...a_{n})$,
puisque cette derni\`{e}re formule est \'{e}quivalent \`{a} $\exists
x_{1}...\exists x_{n}A^{\prime }(x_{1}...x_{n},a_{1}...a_{n}).$ Si on
remplace les $a_{i}$ par des $t_{i}^{\prime }$ avec $t_{i}^{\prime }\equiv
t_{i}[f_{1}:=e_{1}...f_{n}:=e_{n}]$, on obtient une preuve de $\sim
A^{\prime }(e_{1}...e_{n}(t_{1}^{\prime }...t_{n-1}^{\prime }),t_{1}^{\prime
}...t_{n}^{\prime }).$ En appliquant la m\'{e}thode de substitution (qui
reste valable quelles que soient les fonctions calculables avec lesquelles
on a enrichi le langage - maintenant on ne se pr\'{e}occupe plus de
bornage), on obtient bien des fonctions $g_{1}...g_{n}$ \`{a} substituer aux
$e_{1}...e_{n}$ telles que le r\'{e}sultat de $\sim A^{\prime
}(e_{1}...e_{n}(t_{1}^{\prime }...t_{n-1}^{\prime }),t_{1}^{\prime
}...t_{n}^{\prime })$ soit une formule vraie de l'arithm\'{e}tique.

Pourtant ceci ne nous donne pas le r\'{e}sultat recherch\'{e}; en effet, le
fait que $t_{n}$ repr\'{e}sente $\eta _{n}$ ne garantit absolument pas que
l'on ait $t_{n}(e_{1}..e_{n})=\eta _{n}(g_{1}...g_{n}).$ Par exemple, le
terme $\epsilon _{x}(x=y+3)$ est une fonction de $y$ qui repr\'{e}sente la
fonction r\'{e}cursive $y+3$, mais si $\epsilon _{x}(x=y+3)$ appara\^{i}t
dans une preuve, tout ce que nous dit le th\'{e}or\`{e}me 15, c'est que
toutes les formules dans lequel il appara\^{i}t donnent des formules vraies,
mais il se peut tr\`{e}s bien que ce soit le cas sans que pour autant le
substituant final de $\epsilon _{x}(x=y+3)$ soit la fonction r\'{e}cursive $%
\mathbf{+3}$, et on peut m\^{e}me affirmer que ce ne sera pas le cas, car
les fonctions qui jouent le r\^{o}le de substituants sont toujours nulles
sauf pour un nombre fini de valeur. Par contre, comme $\epsilon _{x}(x=y+3)$
repr\'{e}sente $\mathbf{+3}$, il est possible de forcer le r\'{e}sultat pour
certaines valeurs. Il suffit d'ajouter \`{a} l'ensemble de formules de la
preuve des formules comme $\epsilon _{x}(x=t+3)=\mathbf{+3}(t)$ (qu'on
d\'{e}rive de $\forall z(\epsilon _{x}(x=y+3)=\mathbf{+3}(y))$ ce qui force $%
\epsilon _{x}(x=y+3)$ \`{a} se comporter comme $\mathbf{+3}$ pour la valeur
finale de $t$. Ceci est le principe de la solution \`{a} notre probl\`{e}me,
m\^{e}me si la situation est un peu compliqu\'{e}e par le fait qu'on ait
affaire \`{a} des termes repr\'{e}sentant des fonctionnelles et pas
simplement des fonctions.

On se donne deux lemmes qu'on d\'{e}montrera plus tard afin de montrer le
th\'{e}or\`{e}me voulu.

\begin{lemma}
\label{12}Soit $\eta \lbrack a_{1}...a_{n},f_{1}...f_{n}]$ une fonctionnelle
r\'{e}cursive primitive d'ordre fini repr\'{e}sent\'{e}e par un terme $%
t[a_{1}...a_{n},f_{1}...f_{n}]$ et $u_{i}[x_{1}...x_{m_{i}}]$ un terme
quelconque substituable \`{a} $f_{i}$ d'arit\'{e} $m_{i}$, il existe un
terme $R[a_{1}...a_{n}]$ tel que pour tout terme $v[x_{1}...x_{m_{i}}],$ on
peut prouver dans PA$_{\epsilon }$ :

$\left[ \overset{n}{\underset{i=0}{\bigwedge }}\forall \overrightarrow{x}(%
\overrightarrow{x}\leq R[a_{1}...a_{n}]\rightarrow u_{i}[\overrightarrow{x}%
]=v_{i}[\overrightarrow{x}])\right] \rightarrow
t[a_{1}...a_{n},u_{1}...u_{n}]=t[a_{1}...a_{n},v_{1}...v_{n}]$
\end{lemma}

\begin{lemma}
\label{11}Si $\varphi (n)$ est une fonction primitive r\'{e}cursive d'ordre
fini, on peut trouver un terme $t(n)$ tel que l'on peut prouver dans PA$%
_{\epsilon }$ que $t(n)$ satisfait les relations r\'{e}cursives qui
d\'{e}finissent $\varphi (n).$
\end{lemma}

\begin{lemma}
\label{10}On peut trouver des fonctions r\'{e}cursives $h_{1}^{\prime
}...h_{n}^{\prime }$\ param\'{e}tr\'{e}es par des termes $%
u_{i}^{1},v_{i}^{1} $ sans variables libres et des termes $%
u_{i}^{j},v_{i}^{j}$ contenant $N$ comme seule variable libre tels que les
\'{e}galit\'{e}s suivantes soient prouvables dans PA$_{\epsilon }$ :

$e_{1}=h_{1}^{\prime }(u_{1}^{1}...u_{n}^{1},v_{1}^{1}...v_{n}^{1})$

$(y_{1}\leq N\wedge ...\wedge y_{k}\leq N)\rightarrow
e_{k+1}=h_{k+1}^{\prime
}(y_{1}...y_{k},u_{1}^{1}...u_{n}^{k+1},v_{1}^{1}...v_{n}^{k+1})$
\end{lemma}

On a maintenant les outils pour d\'{e}montrer

\begin{theorem}
Si on a une preuve de $\sim A$ dans PA$_{\epsilon },$ alors pour toutes
fonctionnelles primitives r\'{e}cursives d'ordre fini $\eta _{1}...\eta _{n}$
dont les variables sont parmi $f_{1}...f_{n}$, on trouve des fonctions
r\'{e}cursives $g_{1}...g_{n}$ telles que

$\sim A^{\prime }(g_{1}...g_{n}(\overline{\eta }_{1}...\overline{\eta }%
_{n-1}),\overline{\eta }_{1}...\overline{\eta }_{n})$ o\`{u} $\overline{\eta
}_{i}=\eta _{i}[f_{1}:=g_{1},...,f_{n}:=g_{n}]$
\end{theorem}

La d\'{e}monstration commence comme indiqu\'{e} plus haut. Comme par le
lemme \ref{11}, les fonctions r\'{e}cursives primitives d'ordre fini sont
repr\'{e}sentables, il suit directement de la d\'{e}finition des
fonctionnelles primitives r\'{e}cursives d'ordre fini que celles-ci aussi
sont repr\'{e}sentables. Soient donc $t_{1}...t_{n}$ les termes qui les
repr\'{e}sentent. A partir d'une preuve de $\sim A,$ on obtient une preuve de

$\sim A^{\prime }(e_{1}...e_{n}(t_{1}^{\prime }...t_{n-1}^{\prime
}),t_{1}^{\prime }...t_{n}^{\prime })$ \ (0)

o\`{u} les $e_{i}$ et les $t_{i}^{\prime }$ sont d\'{e}finis comme
pr\'{e}c\'{e}demment.

On consid\`{e}re les $h_{i}^{\prime }$ du lemme \ref{10} et on pose $\eta
_{i}^{\prime }=\eta _{i}[f_{1}:=g_{1},...,f_{n}:=g_{n}]$.\ Le lemme \ref{12}
nous dit que si les $\epsilon $-termes $e_{i}$ et les fonctions
r\'{e}cursives $h_{i}^{\prime }$ sont suffisamment \'{e}gaux (au sens de $R$
ind\'{e}pendant des param\`{e}tres de $h_{i}^{\prime }$), les termes
repr\'{e}sentant les fonctionnelles les confondent, ce qui implique que les $%
t_{i}^{\prime }$ et les $\eta _{i}^{\prime }$ donnent \`{a} leur tour le
m\^{e}me r\'{e}sultat (puisque le terme $t_{i}$ auquel on a donn\'{e} comme
argument les $h_{i}^{\prime }$ repr\'{e}sente $\eta _{i}^{\prime }$).
Formellement, le lemme \ref{12} nous donne un terme $R$ tel que

$\left[ \overset{n}{\underset{i=0}{\bigwedge }}\forall \overrightarrow{x}(%
\overrightarrow{x}\leq R\rightarrow e_{i}[\overrightarrow{x}]=h_{i}^{\prime
}[\overrightarrow{x}])\right] \rightarrow \left[ \overset{n}{\underset{i=0}{%
\bigwedge }}t_{i}^{\prime }=\eta _{i}^{\prime }\right] $ \ (1)

Les $t_{i}^{\prime }$ ne comportent pas de variable d'individu libre : la
seule variable libre d'individu $f_{1}$ est remplac\'{e} par le terme clos $%
e_{1}$; par cons\'{e}quent, vu la mani\`{e}re dont $R$ est construit, il ne
comporte pas de variable libre. Et les $t_{i}^{\prime }$ n'ont plus de
variable de fonction libre. Cette remarque est importante car la m\'{e}thode
de substitution s'applique seulement \`{a} des formules closes.
L'introduction de variables de fonctions libres dans le syst\`{e}me sert
\`{a} la formulation des r\'{e}sultats, mais les preuves qui servent \`{a}
leur \'{e}tablissement se font elles sans recourir \`{a} de telles variables.

On pose alors $R^{\prime }=\overset{n}{\underset{i=1}{max}}(R,t_{i}^{\prime
})$

Si on remplace $N$ par $R^{\prime }$ dans le lemme \ref{10}, on des preuves
dans PA$_{\epsilon }$ des formules

$(y_{1}\leq R^{\prime }\wedge ...\wedge y_{k}\leq R^{\prime })\rightarrow
e_{k+1}(y_{1}...y_{k})=$

$h_{k+1}^{\prime
}(y_{1}...y_{k},u_{1}^{1}...u_{n}^{k+1},v_{1}^{1}...v_{n}^{k+1})$ \ (2)

Comme $R_{i}\leq R,$ (2) nous donne les membres de la conjonction qui est
l'ant\'{e}c\'{e}dent de (1) d'o\`{u}

$\overset{n}{\underset{i=0}{\bigwedge }}t_{i}^{\prime }=\eta _{i}^{\prime }$
\ (3)

Dans (3) les $\eta _{i}^{\prime }$ contiennent les $u_{i}^{1},v_{i}^{1}$ et
les $u_{i}^{m+1},v_{i}^{m+1}$ comme param\`{e}tres

Comme $t_{i}^{\prime }\leq R^{\prime },$ (3) nous donne les \'{e}galit\'{e}s

$e_{k+1}(t_{1}^{\prime }...t_{k}^{\prime })=h_{k+1}^{\prime }(\eta
_{1}^{\prime }...\eta _{k}^{\prime
},u_{1}^{1}...u_{n}^{k+1},v_{1}^{1}...v_{n}^{k+1})$ \ (4)

On applique alors la m\'{e}thode de substitution \`{a} l'ensemble de
formules constitu\'{e}es par les preuves de (0), (3) et (4). Le r\'{e}sultat
de (0) en particulier est une formule vraie. Mais on sait en plus que les
substituants des $e_{i}$ co\"{i}ncident avec des fonctions r\'{e}cursives
pour les valeurs des fonctionnelles, et que la valeur pour ces $e_{i}$ des
termes repr\'{e}sentant les fonctionnelles est \'{e}gale \`{a} la valeur des
fonctionnelles pour ces fonctions r\'{e}cursives. Par cons\'{e}quent, si on
appelle $g_{1}...g_{n}$ les substituants des $e_{1}...e_{n},$ (0) nous donne
bien :

$\sim A^{\prime }(g_{1}...g_{n}(\overline{\eta }_{1}...\overline{\eta }%
_{n-1}),\overline{\eta }_{1}...\overline{\eta }_{n})$

\begin{theorem}
Soient deux formules $A\equiv \forall x_{1}\exists y_{1}...\forall
x_{n}\exists y_{n}A^{\prime }(x_{1}...x_{n},y_{1}...y_{n})$

et $B\equiv \forall x_{1}\exists y_{1}...\forall x_{m}\exists y_{m}B^{\prime
}(x_{1}...x_{m},y_{1}...y_{m}),$

une preuve $\pi $\ dans PA$_{\epsilon }$ de $A\rightarrow B$

et $\eta _{1}...\eta _{n}$ des fonctionnelles dont les variables libres sont
parmi $f_{1}...f_{n}$ telles que $A^{\prime }(f_{1}...f_{n}(\eta _{1}...\eta
_{n-1}),\eta _{1}...\eta _{n})$ est v\'{e}rifiable,

on peut trouver des fonctionnelles $\xi _{1}...\xi _{m}$ dont les variables
libres sont parmi $g_{1}...g_{m}$ telles que $B^{\prime }(g_{1}...g_{m}(\xi
_{1}...\xi _{m-1}),\xi _{1}...\xi _{m})$ est v\'{e}rifiable.
\end{theorem}

La d\'{e}monstration suit de pr\`{e}s celle du th\'{e}or\`{e}me
pr\'{e}c\'{e}dent, on garde les notations utilis\'{e}es pour les termes
convoqu\'{e}s dans sa d\'{e}monstration.

Soient $\ \psi _{1}...\psi _{m}$ des fonctions r\'{e}cursives susceptibles
d'instancier les $g_{1}...g_{m}.$

D'une part, on a vu que $\sim A\rightarrow \sim A^{\prime
}(e_{1}...e_{n}(t_{1}^{\prime }...t_{n-1}^{\prime }),t_{1}^{\prime
}...t_{n}^{\prime })$ \'{e}tait prouvable dans PA$_{\epsilon },$ donc par
contrapposition, on a aussi une preuve de $A^{\prime
}(e_{1}...e_{n}(t_{1}^{\prime }...t_{n-1}^{\prime }),t_{1}^{\prime
}...t_{n}^{\prime })\rightarrow A$. D'autre part, comme dans la
d\'{e}monstration du th\'{e}or\`{e}me \ref{100} $B\rightarrow B^{\prime
}(\psi _{1}...\psi _{m}(\epsilon _{1}...\epsilon _{m-1}),\epsilon
_{1}...\epsilon _{m})$ pour certains $\epsilon $-termes $\epsilon
_{1}...\epsilon _{m}.$ Combin\'{e} avec la preuve $\pi $ qui est donn\'{e}e
par l'hypoth\`{e}se, ceci nous donne une preuve de :

$A^{\prime }(e_{1}...e_{n}(t_{1}^{\prime }...t_{n-1}^{\prime
}),t_{1}^{\prime }...t_{n}^{\prime })\rightarrow B^{\prime }(\psi
_{1}...\psi _{m}(\epsilon _{1}...\epsilon _{m-1}),\epsilon _{1}...\epsilon
_{m})$ \ (0)

On peut ajouter aux formules de la preuve de (0) les formules des preuves de
(3) et (4) du th\'{e}or\`{e}me pr\'{e}c\'{e}dent. Ceci nous assure que dans
la substitution finale, le r\'{e}sultat de $A^{\prime
}(e_{1}...e_{n}(t_{1}^{\prime }...t_{n-1}^{\prime }),t_{1}^{\prime
}...t_{n}^{\prime })$ est une formule vraie. Par cons\'{e}quent, le
r\'{e}sultat de $B^{\prime }(\psi _{1}...\psi _{m}(\epsilon _{1}...\epsilon
_{m-1}),\epsilon _{1}...\epsilon _{m})$ est \'{e}galement une formule vraie.
Les valeurs des $\epsilon _{1}...\epsilon _{m}$ sont alors born\'{e}es par $%
born_{m,n,e}^{\prime }[f_{1}:=\psi _{1}...f_{n}:=\psi _{n}]$ et $m,n,e$ ne
d\'{e}pendent pas du choix de $\psi _{1}...\psi _{m}$ (moyennant le choix
d'une preuve standard pour $B\rightarrow B^{\prime }(\psi _{1}...\psi
_{m}(\epsilon _{1}...\epsilon _{m-1}),\epsilon _{1}...\epsilon _{m})$). Il
suit que les $\xi _{i}$ d\'{e}finis par $\xi _{i}=p_{n}^{i}(\epsilon
_{z}(z\leq \beta \wedge B^{\prime
}(g_{1}...g_{m}(p_{m}^{1}(z)...p_{m}^{m}(z)),p_{m}^{1}(z)...p_{m}^{m}(z)))$
o\`{u} $\beta $ est la fonctionnelle $\alpha _{m}(born_{m,n,e}^{\prime
},...,born_{m,n,e}^{\prime })$ conviennent.

\subsection{Retour sur les lemmes}

La d\'{e}monstration du lemme \ref{12} se fait simplement par induction sur
la forme de la fonctionnelle. On indique comment construire $R$ pour chaque
\'{e}tape :

a) max (a$_{i})$

b) max (a$_{i},\chi ,R_{\chi }$)

c) max ($t,R\chi $) o\`{u} $t$ repr\'{e}sente $\omega _{m,\overrightarrow{x}%
}^{\prime }(\chi )$ appliqu\'{e}e aux $u_{i}.$

d) max ($a_{i},R_{\chi _{j}}$)

e) max ($R\chi $)

La d\'{e}monstration du lemme \ref{11} se trouve dans Kreisel \cite{kreis2}.
Intuitivement, la validit\'{e} du lemme est claire : les valeurs des
fonctions primitives r\'{e}cursives d'ordre fini sont calculables en un
nombre fini d'etapes; si toutes les fonctions calculables sont
r\'{e}cursives (th\`{e}se de Church), comme toutes les fonctions
r\'{e}cursives sont repr\'{e}sentables dans AP$_{\epsilon }$ (voir \cite
{hilb1} pour une d\'{e}monstration de ceci) celles-ci le sont \'{e}galement.
On remarque que la notion de repr\'{e}sentation ici utilis\'{e}e est
sp\'{e}cifique \`{a} l'arithm\'{e}tique epsilon : les fonctions sont
repr\'{e}sent\'{e}es par des $\epsilon $-termes et pas par des formules.

Nous donnons enfin la d\'{e}monstration du lemme \ref{10}, un peu difficile
\`{a} suivre quoiqu'assez simple, en suivant scrupuleusement celle de
Kreisel \cite{kreis1}.

On peut trouver des fonctions r\'{e}cursives $h_{1}^{\prime
}...h_{n}^{\prime }$\ param\'{e}tr\'{e}es par des termes $%
u_{i}^{1},v_{i}^{1} $ sans variables libres et des termes $%
u_{i}^{m+1},v_{i}^{m+1}$ contenant $N$ comme seule variable libre tels que

$e_{1}=h_{1}^{\prime }(u_{1}^{1}...u_{n}^{1},v_{1}^{1}...v_{n}^{1})$

$(y_{1}\leq N\wedge ...\wedge y_{k}\leq N)\rightarrow
e_{k+1}=h_{k+1}^{\prime
}(y_{1}...y_{k},u_{k+1}^{k+1}...u_{n}^{k+1},v_{k+1}^{k+1}...v_{n}^{k+1})$

On commence par d\'{e}finir les termes primitifs r\'{e}cursifs (toutes les
minimisations et les quantifications sont born\'{e}es) en jeu. Les $%
a_{j}^{i},b_{j}^{i}$ sont des variables libres

$h_{1}^{\prime }(a_{1}^{1}...a_{n}^{1},b_{1}^{1}...b_{n}^{1})=\epsilon
_{x_{1}}\forall y_{1}...\exists x_{n}\forall y_{n}$

$\left[ y_{1}\leq b_{1}^{1}\wedge ...\wedge y_{n}\leq b_{n}^{1}\rightarrow
x_{1}\leq a_{1}^{1}\wedge ...\wedge x_{n}\leq a_{n}^{1}\sim A^{\prime
}(x_{1}...x_{n},y_{1}...y_{n})\right] $

$%
h_{r+1}(y_{1}...y_{r},c_{1}...c_{r},a_{r+1}^{r+1}...a_{n}^{r+1},b_{r+1}^{r+1}...b_{n}^{r+1})=\epsilon _{x_{r+1}}\forall y_{r+1}...\exists x_{n}\forall y_{n}
$

$\left[ y_{r+1}\leq b_{r+1}^{r+1}\wedge ...\wedge y_{n}\leq
b_{n}^{r+1}\rightarrow x_{r+1}\leq a_{r+1}^{r+1}\wedge ...\wedge x_{n}\leq
a_{n}^{r+1}\sim A^{\prime }(c_{1}...c_{r}x_{r+1}...x_{n},y_{1}...y_{n})%
\right] $

$h_{r+1}^{\prime
}(y_{1}...y_{r},a_{1}^{1}...a_{n}^{r+1},b_{1}^{1}...b_{n}^{r+1})=h_{r+1}(y_{1}...y_{r},h_{1}^{\prime }...h_{r}^{\prime },a_{r+1}^{r+1}...a_{n}^{r+1},b_{r+1}^{r+1}...b_{n}^{r+1})
$

Afin de surmonter la complexit\'{e} de la syntaxe et des doubles indices, on
commence par donner la m\'{e}thode pour trouver les $u_{j}^{i},v_{j}^{i}$
\`{a} partir d'un exemple. On consid\`{e}re la formule $\exists x\forall
yB(a,x,y).$ $a$ est une variable libre qui sert de param\`{e}tre et $B$ peut
comporter d'autres variables li\'{e}es. On consid\`{e}re les termes~:

$\epsilon _{z}\forall v\left[ v\leq N\rightarrow \epsilon _{x}\left( \forall
yB(v,x,y)\right) \leq z\right] $ qu'on note $u$

$\epsilon _{z}\forall v,w\left[ v\leq N\wedge w\leq u\rightarrow \epsilon
_{y}\left( \sim B(v,w,y)\right) \leq z\right] $ qu'on note $v$

Intuitivement $u$ et $v$ sont deux des $u_{j}^{i},v_{j}^{i}$ qu'on cherche,
c'est-\`{a}-dire des param\`{e}tres qui, relativement \`{a} un $N$
quelconque permettent de borner les quantifications. C'est ce que veulent
dire les deux propositions suivantes :

\begin{proposition}
\label{335}$a\leq N\rightarrow $ $\epsilon _{x}\left( \forall
yB(a,x,y)\right) =\epsilon _{x}\forall y\left[ y\leq v\rightarrow x\leq
u\wedge B(a,x,y)\right] $
\end{proposition}

\begin{proposition}
\label{336}$a\leq N\rightarrow $ $\exists x\forall
yB(a,x,y)\longleftrightarrow \exists x\forall y\left[ y\leq v\rightarrow
x\leq u\wedge B(a,x,y)\right] $
\end{proposition}

En toute rigueur, on aurait besoin pour le lemme d'\'{e}tablir que ces
propositions sont prouvables dans PA$_{\epsilon }.$ On se contente ici de
les prouver informellement, la formalisation, fastidieuse, ne posant pas de
probl\`{e}me particulier.

$\epsilon _{x}\left( \forall yB(a,x,y)\right) $ est not\'{e} $\lambda (a)$

$\epsilon _{x}\forall y\left[ y\leq v\rightarrow x\leq u\wedge B(a,x,y)%
\right] $ est not\'{e} $\sigma (a)$

$\epsilon _{y}\sim B(a,c,y)$ est not\'{e} $\rho (a,c)$

Premier cas : $\exists x\forall yB(a,x,y)$ est vrai pour tout $a\leq N$. $%
\forall yB(a,\lambda (a),y)$ est donc vrai aussi. Par ailleurs, quand $a$
parcourt [0;N], il existe un maximum $x$ des $\lambda (a)$ et $x\leq u.$
Donc \ref{336} est vrai. Mais est-ce que $\epsilon _{x}\forall y\left[ y\leq
v\rightarrow x\leq u\wedge B(a,x,y)\right] $ pourrait \^{e}tre strictement
inf\'{e}rieur \`{a} cet $x$ ? Soit un $\ c$ qui pour un certain $a$\ est
strictement inf\'{e}rieur \`{a}$\ \lambda (a),$ il y a donc un $y\in \{\rho
(a,c)/a\leq N,c\leq \lambda (a)\}$ tel que $\sim B(a,c,y),$
c'est--\`{a}-dire un $y\leq v$ puisque $v$ majore cet ensemble. Autrement
dit, $c$ ne peut \^{e}tre $\epsilon _{x}\forall y\left[ y\leq v\rightarrow
x\leq u\wedge B(a,x,y)\right] .$ On a ainsi montr\'{e} \ref{335}.

Deuxi\`{e}me cas : il existe un $a\leq N$ tel que $\forall x\exists
yB(a,x,y) $. On a ainsi $\forall x\sim B(a,x,\rho (a,x)).$ Pour $x\leq u,$ $%
\rho (a,x)\leq v$ par d\'{e}finition de $v$. Il n'y a donc pas de $x\leq u$
tel que $\forall y\leq v\left( B(a,x,y)\right) .$ Les deux membres de
l'\'{e}quivalence de \ref{336} sont donc faux et les deux $\epsilon $-termes
de \ref{335} nuls.

On montre maintenant comment utiliser la m\'{e}thode de l'exemple pour
obtenir les $u_{i}^{1},v_{i}^{1}.$

On commence par appliquer la m\'{e}thode\ avec $B(a,x,y)\equiv \exists
x_{2}\forall y_{2}...\exists x_{n}\forall y_{n}\sim A^{\prime
}(x,x_{2}...x_{n},y,y_{2}...y_{n})$ pour avoir $u_{1}^{1}\equiv
u,v_{1}^{1}\equiv v$. La proposition \ref{335} nous donne alors

$\epsilon _{x_{1}}\left( \forall y_{1}...\exists x_{n}\forall y_{n}\sim
A^{\prime }(x_{1}...x_{n},y_{1}...y_{n})\right) $

$=\epsilon _{x_{1}}\forall y_{1}\left[ y_{1}\leq v_{1}^{1}\rightarrow
x_{1}\leq u_{1}^{1}\wedge \exists x_{2}\forall y_{2}...\exists x_{n}\forall
y_{n}\sim A^{\prime }(x_{1},...x_{n},y_{1}...y_{n})\right] $

Et on continue : on proc\`{e}de de m\^{e}me avec

$B(x_{1},y_{1},x_{2},y_{2})\equiv \exists x_{2}\forall y_{2}\left[ y_{1}\leq
v_{1}^{1}\rightarrow x_{1}\leq u_{1}^{1}\wedge \exists x_{3}\forall
y_{3}...\exists x_{n}\forall y_{n}\sim A^{\prime
}(x_{1},...x_{n},y_{1}...y_{n})\right] $

pour obtenir $u_{2}^{1}\equiv u,v_{2}^{1}\equiv v$. Il faut pr\'{e}ciser les
$N$ qui bornent les $x_{1},y_{1}$ qui jouent ici le r\^{o}le des
param\`{e}tres $a.$ On exige $x_{1}\leq u_{1}^{1}$ et $y_{1}\leq v_{1}^{1}.$

On doit alors prouver

$\epsilon _{x_{1}}\left( \forall y_{1}...\exists x_{n}\forall y_{n}\sim
A^{\prime }(x_{1}...x_{n},y_{1}...y_{n})\right) $

$=\epsilon _{x_{1}}\forall y_{1}\exists x_{2}\forall _{2}\left[ y_{1}\leq
v_{1}^{1}\wedge y_{2}\leq v_{2}^{1}\rightarrow x_{1}\leq u_{1}^{1}\wedge
x_{2}\leq u_{2}^{1}\wedge \exists x_{2}\forall y_{2}...\exists x_{n}\forall
y_{n}\sim A^{\prime }(x_{1},...x_{n},y_{1}...y_{n})\right] $

Pour cela, il suffit d'avoir

$\epsilon _{x_{1}}\forall y_{1}\left[ y_{1}\leq v_{1}^{1}\rightarrow
x_{1}\leq u_{1}^{1}\wedge \exists x_{2}\forall y_{2}...\exists x_{n}\forall
y_{n}\sim A^{\prime }(x_{1},...x_{n},y_{1}...y_{n})\right] $

$=\epsilon _{x_{1}}\forall y_{1}\exists x_{2}\forall _{2}\left[ y_{1}\leq
v_{1}^{1}\wedge y_{2}\leq v_{2}^{1}\rightarrow x_{1}\leq u_{1}^{1}\wedge
x_{2}\leq u_{2}^{1}\wedge \exists x_{2}\forall y_{2}...\exists x_{n}\forall
y_{n}\sim A^{\prime }(x_{1},...x_{n},y_{1}...y_{n})\right] $

On peut alors se ramener \`{a} montrer dans PA$_{\epsilon }$
l'\'{e}quivalence

$\exists x_{2}\forall y_{2}\left[ y_{1}\leq v_{1}^{1}\rightarrow x_{1}\leq
u_{1}^{1}\wedge \exists x_{3}\forall y_{3}...\exists x_{n}\forall y_{n}\sim
A^{\prime }(x_{1},...x_{n},y_{1}...y_{n})\right] $

$\longleftrightarrow \exists x_{2}\forall y_{2}\left[ y_{1}\leq
v_{1}^{1}\wedge y_{2}\leq v_{2}^{1}\rightarrow x_{1}\leq u_{1}^{1}\wedge
x_{2}\leq u_{2}^{1}\wedge \exists x_{2}\forall y_{2}...\exists x_{n}\forall
y_{n}\sim A^{\prime }(x_{1},...x_{n},y_{1}...y_{n})\right] $

On raisonne par cas : si $y_{1}>v_{1}^{1},$ les deux formules sont vraies.
Si $y_{1}\leq v_{1}^{1}$ et $x_{1}>u_{1}^{1}$ la premi\`{e}re formule est
fausse, et la seconde aussi (prendre $y_{2}=0$ pour s'assurer que
l'ant\'{e}c\'{e}dent est vrai). Si $y_{1}\leq v_{1}^{1}$ et $x_{1}\leq
u_{1}^{1},$ l'ant\'{e}c\'{e}dent du lemme \ref{336} est vrai ce qui nous
donne

$\exists x_{2}\forall y_{2}\exists x_{3}\forall y_{3}...\exists x_{n}\forall
y_{n}\sim A^{\prime }(x_{1},...x_{n},y_{1}...y_{n})$

$\longleftrightarrow \exists x_{2}\forall y_{2}\left[ y_{2}\leq
v_{2}^{1}\rightarrow x_{2}\leq u_{2}^{1}\wedge \exists x_{2}\forall
y_{2}...\exists x_{n}\forall y_{n}\sim A^{\prime
}(x_{1},...x_{n},y_{1}...y_{n})\right] $

et l'\'{e}quivalence du m\^{e}me coup.

On obtient les autres $u_{i}^{1},v_{i}^{1}$ en r\'{e}p\'{e}tant la
proc\'{e}dure. On a ainsi un $h_{1}^{\prime
}(u_{1}^{1}...u_{n}^{1},v_{1}^{1}...v_{n}^{1})$ prouvablement \'{e}gal \`{a}
$e_{1}.$

On montre maintenant comment obtenir $h_{r+1}^{\prime
}(y_{1}...y_{r},u_{1}^{1}...u_{n}^{r+1},v_{1}^{1}...v_{n}^{r+1})$ tel que

$(y_{1}\leq N\wedge ...\wedge y_{k}\leq N)\rightarrow
e_{r+1}=h_{r+1}^{\prime
}(y_{1}...y_{r},u_{1}^{1}...u_{n}^{r+1},v_{1}^{1}...v_{n}^{r+1})$

en supposant qu'on a su le faire pour $h_{1}^{\prime }...h_{r}^{\prime }.$

On consid\`{e}re le terme $t\equiv \epsilon _{x_{r+1}}\forall
y_{r+1}...\exists x_{n}\forall y_{n}\sim A^{\prime
}(c_{1}...c_{r}x_{r+1}...x_{n},y_{1}...y_{n})$ dont les variables sont $%
c_{1}...c_{r},y_{1}...y_{r}$

On sait trouver des termes $u_{r+1}...u_{n},$ $v_{r+1}...v_{r}$ contenant
les variables libres $c_{1}...c_{r}$ et $N$\ tels que l'on peut prouver :

$(y_{1}\leq N\wedge ...\wedge y_{k}\leq N)\rightarrow t=\epsilon
_{x_{r+1}}\forall y_{r+1}...\exists x_{n}\forall y_{n}$

$\left[ y_{r+1}\leq v_{r+1}\wedge ...\wedge y_{n}\leq v_{n}\rightarrow
x_{r+1}\leq u_{r+1}\wedge ...\wedge x_{n}\leq a_{n}\sim A^{\prime
}(c_{1}...c_{r},x_{r+1}...x_{n},y_{1}...y_{n})\right] $

On pose alors, pour $r+1\leq j\leq n$, $u_{j}^{r+1}=u_{j}\left[
c_{1}:=h_{1}^{\prime }...c_{r}:=h_{r}^{\prime }\right] $ et de m\^{e}me $%
v_{j}^{r+1}=v_{j}\left[ c_{1}:=h_{1}^{\prime }...c_{r}:=h_{r}^{\prime }%
\right] .$

Dans l'\'{e}galit\'{e} pr\'{e}c\'{e}dente, on remplace dans $t$ les $c_{j}$
par les $e_{j}$ ce qui nous donne $e_{r+1},$ et dans le membre de droite les
$c_{j}$ par les $h_{j}^{\prime }.$ Cette substitution est licite car pour $%
y_{1}\leq N\wedge ...\wedge y_{k}\leq N,$ on sait que $e_{j}=h_{j}^{\prime
}. $ Ceci nous donne alors~:

$(y_{1}\leq N\wedge ...\wedge y_{k}\leq N)\rightarrow
e_{r+1}=h_{r+1}^{\prime
}(y_{1}...y_{r},u_{1}^{1}...u_{n}^{r+1},v_{1}^{1}...v_{n}^{r+1}).$

\section{R\'{e}alisabilit\'{e} classique et arithm\'{e}tique}

\subsection{Le lambda-c calcul et les r\`{e}gles de typage}

On d\'{e}finit l'ensemble $\Lambda _{c}$ des termes t et l'ensemble $\Pi $
des piles $\pi $ :

$t=cc,x,(t)t,\lambda x.t,k_{\pi }$

$\pi =\rho ,t.\pi $

La r\'{e}duction $\succ $\ est alors d\'{e}finie sur $\Lambda _{c}\times \Pi
$, l'ensemble des ex\'{e}cutables~:

$(t)u\star \pi \succ t\star u\cdot \pi \qquad \qquad \ $(push)

$\lambda x.t\star u\cdot \pi \succ t[u/x]\star \pi \qquad $(pop)

$cc\star t\cdot \pi \succ t\star k_{\pi }\cdot \pi \qquad $\ \ \ \ (store)

$k_{\pi }\star t\cdot \pi ^{\prime }\succ t\star \pi \qquad \qquad \ \ \ $%
(restore)

Les r\`{e}gles de typage correspondent \`{a} la logique classique du second
ordre. Les seuls symboles logiques utilis\'{e}s sont $\rightarrow $ et $%
\forall $.

Soit $\Gamma $ un contexte de la forme $x_{1}:A_{1},...,x_{n}:A_{n}$, les
r\`{e}gles de typage sont~:

1. $\Gamma \vdash x_{i}:A_{i}$ \ o\`{u} $1\leq i\leq n$

2. si $\Gamma \vdash t:A\rightarrow B$ et $\Gamma \vdash u:A$ alors $\Gamma
\vdash tu:B$

3. si $\Gamma ,x:A\vdash t:B$ alors $\Gamma \vdash \lambda x.t:A\rightarrow
B $

4. si $\Gamma \vdash t:(A\rightarrow B)\rightarrow A$ alors $\Gamma \vdash
cct:A$

5. si $\Gamma \vdash t:A$ alors $\Gamma \vdash t:\forall xA$ pour $x$ non
libre dans $\Gamma $

6. si $\Gamma \vdash t:A$ alors $\Gamma \vdash t:\forall XA$ pour $X$ non
libre dans $\Gamma $

7. si $\Gamma \vdash t:\forall xA$ alors $\Gamma \vdash t:A[\tau /x]$ pour
tout terme $\tau $

8. si $\Gamma \vdash t:\forall XA$ alors $\Gamma \vdash t:A[\phi
(x_{1}...x_{n})/Xx_{1}...x_{n}]$ pour toute formule $\phi (x_{1}...x_{n})$
(axiome de compr\'{e}hension pour la logique du second ordre)

Les autres connecteurs sont d\'{e}finis ensuite de la mani\`{e}re usuelle :

$\perp \equiv \forall XX$

$A\wedge B\equiv \forall X((A\rightarrow (B\rightarrow X))\rightarrow X)$

$A\vee B\equiv \forall X(((A\rightarrow X)\rightarrow (B\rightarrow
X))\rightarrow X)$

$\exists XA\equiv \forall X(A\rightarrow \perp )\rightarrow \perp $

De plus $a=b\equiv \forall X(Xa\rightarrow Xb)$

\subsection{La r\'{e}alisabilit\'{e}}

L'id\'{e}e de la r\'{e}alisabilit\'{e} est d'interpr\'{e}ter les formules F
par des sous-ensembles $\left| F\right| $ de $\Lambda _{c}.$ La
d\'{e}termination de $\left| F\right| $ se fait relativement \`{a} un
mod\`{e}le M de la logique classique du second ordre et \`{a} un ensemble $%
\amalg \sqsubseteq \Lambda _{c}\times \Pi $ de processus qui correspond aux
processus observables. On montre des propri\'{e}t\'{e}s sur les termes qui
sont dans $\left| F\right| $ (qui r\'{e}alisent F) en choisissant des $%
\amalg $ particuliers.\ Ceci permet ensuite d'avoir ces propri\'{e}t\'{e}s
pour les termes de type F, c'est-\`{a}-dire pour les termes qui codent des
preuves de F moyennant un lemme d'ad\'{e}quation qui nous assure que si $%
\vdash t:A$ alors $t\in \left| A\right| .$

Un mod\`{e}le $M$ est d\'{e}fini comme la donn\'{e}e d'un ensemble M
d'individus et d'une interpr\'{e}tation $f_{M}:M^{k}\mapsto M$ pour chaque
symbole de fonction f k-aire du langage. Le domaine de variation des
pr\'{e}dicats d'arit\'{e} k est $P(\Pi )^{M^{k}}.$ On exige que $\amalg $
soit clos par anti-r\'{e}duction.

On d\'{e}finit la valeur de v\'{e}rit\'{e} d'une formule F par

$\left| A\right| =\left\| A\right\| ^{\perp }=\left\{ t\in \Lambda
_{c}/\forall \pi _{A}\in \left\| A\right\| ,t\ast \pi _{A}\in \amalg
\right\} $

On d\'{e}finit maintenant $\left\| A\right\| \subseteq \Pi $ pour une
formule $A$ \`{a} param\`{e}tre dans M par induction sur $A.$

Si $A$ est atomique $A\equiv Ra_{1}...a_{k}$ et $R\in P(\Pi )^{M^{k}},$ on
pose $\left\| Ra_{1}...a_{k}\right\| =R(a_{1}...a_{k})$

$\left\| A\rightarrow B\right\| =\{t.\pi /t\in \left| A\right| ,\pi \in
\left\| B\right\| \}$

$\left\| \forall xA\right\| =\underset{a\in M}{\cup }\left\| A[x:=a]\right\|
$

$\left\| \forall X^{k}A\right\| =\underset{R\in P(\Pi )^{M^{k}}}{\cup }%
\left\| A[X:=R]\right\| $

Quelques valeurs sont distingu\'{e}es. $\top =\parallel \emptyset \parallel
^{\perp }$ est la plus grande valeur de v\'{e}rit\'{e}. $\perp =\parallel
\Pi \parallel ,$ donc $\left| \perp \right| $ la plus petite valeur de
v\'{e}rit\'{e}.

\begin{theorem}
Si $x_{1}:F_{1}...x_{n}:F_{n}\vdash t:F$ et si $s_{1}\in $ $\left|
F_{1}\right| ,...,s_{n}\in $ $\left| F_{n}\right| $ alors $%
t[s_{1}/x_{1}...s_{n}/x_{n}]\in \left| F\right| $
\end{theorem}

La preuve se fait par induction sur le typage. On examine les diff\'{e}rents
cas :

1. La derni\`{e}re r\`{e}gle est : $\Gamma \vdash x_{i}:A_{i}$ \ o\`{u} $%
1\leq i\leq n,$ c'est imm\'{e}diat.

2. La derni\`{e}re r\`{e}gle est : si $\Gamma \vdash t:A\rightarrow B$ et $%
\Gamma \vdash u:A$ alors $\Gamma \vdash tu:B$

$tu[\overrightarrow{s}/\overrightarrow{x}]=t[\overrightarrow{s}/%
\overrightarrow{x}]u[\overrightarrow{s}/\overrightarrow{x}]$ et
l'hypoth\`{e}se d'induction nous donne $t[\overrightarrow{s}/\overrightarrow{%
x}]\in A\rightarrow B$ et $u[\overrightarrow{s}/\overrightarrow{x}]\in
\left| A\right| .$ Donc $t[\overrightarrow{s}/\overrightarrow{x}]\ast u[%
\overrightarrow{s}/\overrightarrow{x}].\pi _{B}\in \amalg $ donc
(cl\^{o}ture par anti-r\'{e}duction de $\amalg $), $t[\overrightarrow{s}/%
\overrightarrow{x}]u[\overrightarrow{s}/\overrightarrow{x}]\ast \pi _{B}\in
\amalg ,$ d'o\`{u} $tu[\overrightarrow{s}/\overrightarrow{x}]\in \left|
B\right| $

3. La derni\`{e}re r\`{e}gle est : si $\Gamma ,y:A\vdash t:B$ alors $\Gamma
\vdash \lambda y.t:A\rightarrow B.$

On suppose que $y$ est une variable fraiche qui n'appara\^{i}t pas dans les $%
s_{i}$ et qui est diff\'{e}rente des $x_{i}.$ $(\lambda y.t)[\overrightarrow{%
s}/\overrightarrow{x}]=\lambda y.t[\overrightarrow{s}/\overrightarrow{x}].$
Soient $u\in \left| A\right| ,$ $\pi _{B}\in \left\| B\right\| ,$ on veut $%
\lambda y.t[\overrightarrow{s}/\overrightarrow{x}]\ast u.\pi _{B}\in \amalg
. $

$\lambda y.t[\overrightarrow{s}/\overrightarrow{x}]\ast u.\pi _{B}\succ t[%
\overrightarrow{s}/\overrightarrow{x},u/y]\ast \pi _{B}$ or l'hypoth\`{e}se
de r\'{e}currence nous donne pr\'{e}cis\'{e}ment $[\overrightarrow{s}/%
\overrightarrow{x},u/y]\ast \pi _{B}\in \amalg ,$ d'o\`{u} le r\'{e}sultat.

4. La derni\`{e}re r\`{e}gle est : si $\Gamma \vdash t:(A\rightarrow
B)\rightarrow A$ alors $\Gamma \vdash cct:A$

On veut $cc\in \left| ((A\rightarrow B)\rightarrow A)\rightarrow A\right| .$
Supposons $u_{(A\rightarrow B)\rightarrow A}\in \left| (A\rightarrow
B)\rightarrow A\right| $ et $\pi _{A}\in \left\| A\right\| .$ On veut $%
cc\ast u_{(A\rightarrow B)\rightarrow A}.\pi _{A}\in \amalg .$ $cc\ast
u_{(A\rightarrow B)\rightarrow A}.\pi _{A}\succ u_{(A\rightarrow
B)\rightarrow A}\ast k_{\pi _{A}}.\pi _{A}.$ Ceci est dans $\amalg $ \`{a}
condition que $k_{\pi _{A}}\in \left| A\rightarrow B\right| .$ Soient $%
v_{A}\in \left| A\right| $ et $\pi _{B}\in \left\| B\right\| ,$ est-ce que $%
k_{\pi _{A}}\ast v_{A}.\pi _{B}\in \amalg $ ? oui, car $k_{\pi _{A}}\ast
v_{A}.\pi _{B}\succ v_{A}\ast k_{\pi _{A}}.$

5. La derni\`{e}re r\`{e}gle est : si $\Gamma \vdash t:A$ alors $\Gamma
\vdash t:\forall xA$ pour $x$ non libre dans $\Gamma .$

L'hypoth\`{e}se d'induction nous donne $t[\overrightarrow{s}/\overrightarrow{%
x}]\in \left| A[b/x]\right| $ pour tout $b$ donc $t[\overrightarrow{s}/%
\overrightarrow{x}]\in \underset{b\in M}{\cap }\left| A[b/x]\right| =\left|
\forall xA\right| $

6. La derni\`{e}re r\`{e}gle est : si $\Gamma \vdash t:A$ alors $\Gamma
\vdash t:\forall X^{k}A$ pour $X^{k}$ non libre dans $\Gamma .$

L'hypoth\`{e}se d'induction nous donne $t[\overrightarrow{s}/\overrightarrow{%
x}]\in \left| A[R^{k}/X^{k}]\right| $ pour tout $R^{k}$ donc $t[%
\overrightarrow{s}/\overrightarrow{x}]\in \underset{R\in P(\Pi )^{M^{k}}}{%
\cap }\left| A[R^{k}/X^{k}]\right| =\left| \forall X^{k}A\right| $

7. si $\Gamma \vdash t:\forall xA$ alors $\Gamma \vdash t:A[\tau /x]$ un
terme $\tau $

Il suffit de voir que $\left| A[\tau /x]\right| \sqsupseteq \left| \forall
xA\right| $ (consid\'{e}rer $A[b/x]$ o\`{u} $b=\tau )$ et d'appliquer
l'hypoth\`{e}se d'induction.

8. La derni\`{e}re r\`{e}gle est : si $\Gamma \vdash t:\forall XA$ alors $%
\Gamma \vdash t:A[\phi (x_{1}...x_{n})/Xx_{1}...x_{n}]$ pour une formule $%
\phi (x_{1}...x_{n})$

Il suffit de voir que $\left| A[\phi /X]\right| \sqsupseteq \left| \forall
XA\right| $ (consid\'{e}rer $A[R/\phi ]$ o\`{u} $R(a_{1}...a_{n})=\left\|
\phi (a_{1}...a_{n})\right\| $ et d'appliquer l'hypoth\`{e}se d'induction.

\begin{corollary}
Si $\vdash t:F$ alors $t\in \left| F\right| $
\end{corollary}

On montre \'{e}galement un th\'{e}or\`{e}me simple sur l'\'{e}galit\'{e} qui
sera utile par la suite.

\begin{theorem}
\label{4}Soient a,b deux entiers;

$\Vert a=b\Vert =\Vert \top \rightarrow \bot \Vert $ si $a\neq b.$

$\Vert a=b\Vert =\Vert \forall X\left( X\rightarrow X\right) \Vert $ si $%
a=b. $
\end{theorem}

$\Vert a=b\Vert =_{def}\left\| \forall X(Xa\rightarrow Xb)\right\| $

Si $a\neq b,$ on peut prendre $\left| Xa\right| =\Lambda _{c}$ et $\left\|
Xb\right\| =\Pi $ d'o\`{u} $\Vert a=b\Vert \supseteq \{t.\pi /t\in \Lambda
_{c},\pi \in \Pi \}$ et l'\'{e}galit\'{e} suit trivialement de ce que les
piles de $\Vert a=b\Vert $ ayant n\'{e}cessairement la forme $t.\pi ,$ il
n'existe par de surensemble de piles de cette forme.

Si $a=b,$ on a n\'{e}cessairement $\left\| Xa\right\| =\left\| Xb\right\| ,$
d'o\`{u} $\left\| \forall X(Xa\rightarrow Xb)\right\| =$ $\Vert \forall
X\left( X\rightarrow X\right) \Vert .$

\subsection{L'arithm\'{e}tique du second ordre}

\begin{definition}
Une formule F est r\'{e}alis\'{e}e s'il existe un lambda-c terme t sans
continuation tel que $t\Vdash F$ pour tout choix de $\amalg .$
\end{definition}

On part d'une formulation standard de l'arithm\'{e}tique de Peano du second
ordre (PA$_{2}$), par exemple l'ensemble d'axiomes suivants formul\'{e}s
dans un langage L contenant $\{0,s,+\times \}$ :

1. axiomes pour le successeur

$s0\neq 0$ et $\forall x,y\left( sx=sy\rightarrow x=y\right) $

2. axiomes pour l'addition

$\forall x(x+0=x)$

$\forall x,y(x+sy=s(x+y))$

3. axiomes pour la multiplication

$\forall x(x\times 0=0)$

$\forall x,y(x\times sy=x\times y+y)$

4. axiome d'induction

$\forall xInt(x)$ o\`{u} $Int(x)\equiv \forall X\left[ \forall y\left(
Xy\rightarrow Xsy\right) ,X0\rightarrow Xx\right] $

\bigskip

Afin de pouvoir typer tous les th\'{e}or\`{e}mes de PA$_{2},$ on voudrait
que tous ces axiomes soient r\'{e}alis\'{e}s. En effet, si un axiome $A$ est
r\'{e}alis\'{e}, disons par un terme $t$, alors le typage $\Gamma \vdash t:A$
est compatible avec le lemme d'ad\'{e}quation.

1. Soit $u$ un terme quelconque. On montre que $\lambda x.xu\in \left|
s0\neq 0\right| =\left| (\top \rightarrow \bot )\rightarrow \bot \right| $.
Soit $t_{\top \rightarrow \bot }\in \left| \top \rightarrow \bot \right| ,$
pour toute pile $\pi ,$ $t_{\top \rightarrow \bot }\ast u.\pi \in \amalg .$
Or $\lambda x.xu\ast t_{\top \rightarrow \bot }.\pi \succ $ $t_{\top
\rightarrow \bot }\ast u.\pi $. Par ailleurs, il est clair que $\lambda
x.x\in \left| \forall x,y\left( sx=sy\rightarrow x=y\right) \right| .$

2 et 3. Il suit du th\'{e}or\`{e}me \ref{4} que toutes les formules
\'{e}quationnelles vraies sont r\'{e}alis\'{e}es par l'identit\'{e}.

Par contre, l'axiome d'induction n'est pas r\'{e}alis\'{e} (on montre qu'il
ne l'est pas m\^{e}me pour un mod\`{e}le \`{a} deux \'{e}l\'{e}ments). Ceci
conduit \`{a} envisager une relativisation des quantificateurs du premier
ordre : on cherche \`{a} se dispenser de l'axiome d'induction en limitant
les th\'{e}or\`{e}mes aux \'{e}l\'{e}ments qui le satisfont. On devra
utiliser pour cela :

4'. $\forall x_{1}...\forall x_{k}\left[ Int(x_{1})...Int(x_{k})\rightarrow
Int(fx_{1}...x_{k})\right] $ pour tout symbole de fonction $f$ du langage.

\begin{definition}
On d\'{e}finit $A^{int}$ par induction sur la forme de $A$

si $A$ est atomique, $A^{int}\equiv A$

$(A\rightarrow B)^{int}\equiv A^{int}\rightarrow B^{int}$

$(\forall xA)^{int}\equiv \forall x(Int(x)\rightarrow A^{int})$

$(\forall XA)^{int}\equiv \forall X(A^{int})$
\end{definition}

\begin{theorem}
Si $\vdash _{PA_{2}}A,$ alors $\vdash _{1+2+3+4^{\prime }}A^{int}$ pour $A$
une formule close.
\end{theorem}

L'induction sur la longueur des preuves exige une formulation un peu plus
g\'{e}n\'{e}rale.

\begin{theorem}
Si $\Gamma \vdash _{PA_{2}}A,$ alors $\Gamma ^{int}\cup \left\{
Int(x_{i})\right\} _{x_{i}\in vlib(A)}\vdash _{1+2+3+4^{\prime }}A^{int}.$
\end{theorem}

La d\'{e}monstration par induction ne pose pas de probl\`{e}me particulier.
On s'attarde seulement sur les r\`{e}gles o\`{u} quelque chose se passe.

a) Si $\Gamma \vdash A$ alors $\Gamma \vdash t:\forall xA$ pour $x$ non
libre dans $\Gamma $

L'hypoth\`{e}se d'induction nous dit $\Gamma ^{int}\cup \left\{
Int(x_{i})\right\} _{x_{i}\in vlib(A)}\vdash _{1+2+3+4^{\prime }}A^{int}.$

On veut $\Gamma ^{int}\cup \left\{ Int(x_{i})\right\} _{x_{i}\in
vlib(\forall xA)}\vdash _{1+2+3+4^{\prime }}\forall x(Int(x)\rightarrow
A^{int})$

Il suffit d'appliquer une $\rightarrow _{intro}$ qui \'{e}limine $Int(x)$
des hypoth\`{e}ses.

b) si $\Gamma \vdash \forall xA$ alors $\Gamma \vdash t:A[\tau /x]$ pour
tout terme $\tau .$

L'hypoth\`{e}se d'induction nous donne $\Gamma ^{int}\cup \left\{
Int(x_{i})\right\} _{x_{i}\in vlib(\forall xA)}\vdash _{1+2+3+4^{\prime
}}\forall x(Int(x)\rightarrow A^{int})$

L'id\'{e}e est d'effectuer une $\forall _{\acute{e}lim}$ qui introduise $%
\tau $ puis une $\rightarrow _{\acute{e}lim}$ avec $Int(\tau ).$ On montre
facilement par induction sur la forme d'un terme qu'on peut avoir $Int(\tau
) $ pour n'importe quel terme, en utilisant $Int(0),$ 4' et des
hypoth\`{e}ses suppl\'{e}mentaires $Int(x_{j})$ pour les variables libres $%
x_{j}$ de $\tau .$ Donc on aura bien :

$\Gamma ^{int}\cup \left\{ Int(x_{i})\right\} _{x_{i}\in vlib(A[\tau
/x])}\vdash _{1+2+3+4^{\prime }}A[\tau /x]^{int}).$

c) On v\'{e}rifie enfin facilement que $\vdash \forall xInt(x)^{int}$

Reste alors \`{a} montrer

\begin{theorem}
\label{5}$\forall x_{1}...\forall x_{k}\left[ Int(x_{1})...Int(x_{k})%
\rightarrow Int(fx_{1}...x_{k})\right] $ est r\'{e}alis\'{e} pour tout
symbole de fonction $f$ du langage.
\end{theorem}

\begin{lemma}
\label{6}Soient $\xi ,\eta ,t_{1}...t_{k}$ des $\lambda $-termes du $\lambda
$-calcul ordinaire, si $\eta $ n'est pas une application et si $\xi >\eta
t_{1}...t_{k}$ o\`{u} $>$ d\'{e}signe la r\'{e}duction de t\^{e}te
paresseuse, alors pour toute pile $\pi ,$ $\xi \ast \pi \succ \eta
t_{1}...t_{k}\ast \pi $
\end{lemma}

La d\'{e}monstration se fait par induction sur la longueur de la
r\'{e}duction de t\^{e}te. Si $\xi $ est $\eta t_{1}...t_{k},$ le
r\'{e}sultat est imm\'{e}diat. Si $\xi $ est $(\lambda yu)v\overrightarrow{w}%
,$ $\xi >_{1}u[y:=v]\overrightarrow{w}>\eta t_{1}...t_{k}.$ Par
hypoth\`{e}se d'induction, $u[y:=v]\overrightarrow{w}\ast \pi \succ \eta
t_{1}...t_{k}\ast \pi .$ Comme $\eta $ n'est pas une application, des
\'{e}tapes de r\'{e}duction doivent n\'{e}cessairement avoir lieu qui
empilent les $\overrightarrow{w}$ sur $\pi .$ Donc $u[y:=v]\overrightarrow{w}%
\ast \pi \succ u[y:=v]\ast \overrightarrow{w}.\pi \succ \eta
t_{1}...t_{k}\ast \pi .$ D'o\`{u} $\xi \ast \pi \succ u[y:=v]\ast
\overrightarrow{w}.\pi \succ \eta t_{1}...t_{k}\ast \pi .$

On ne red\'{e}montre pas le th\'{e}or\`{e}me suivant qui donne pour la
d\'{e}finition de la r\'{e}alisabilit\'{e} classique un r\'{e}sultat du $%
\lambda $-calcul ordinaire.

\begin{theorem}
\label{8}Soit $n$ un entier et $\nu $ un $\lambda $-terme du $\lambda $%
-calcul ordinaire $\beta $-\'{e}quivalent \`{a} l'entier de Church $n,$ $\nu
\Vdash Int[s^{n}0]$
\end{theorem}

Soit s un $\lambda $-terme d\'{e}termin\'{e} repr\'{e}sentant le successeur
sur les entiers de Church. On d\'{e}finit $T=\lambda f\lambda n.(((n)\lambda
gg\circ s)f)0$ (T est l'op\'{e}rateur de stockage de \cite{kr3}). On montre
le lemme suivant

\begin{lemma}
\label{9}si $\phi \ast s^{n}0.\pi _{X}$ pour tout $\pi _{X}\in \left\|
X\right\| ,$ alors $T\phi \Vdash Int(n)\rightarrow X$
\end{lemma}

Soit $\nu \Vdash Int(n).$ Soit $\amalg $ un choix de bottom fix\'{e} et $\pi
_{X}$ une pile quelconque dans $\Vert X\parallel ,$ on veut montrer

$T\phi \ast \nu .\pi _{X}\in \amalg .$

$T\phi \ast \nu .\pi _{X}\succ \lambda f\lambda n.(((n)\lambda gg\circ
s)f)0\ast \phi .\nu .\pi _{X}\succ ((\nu )\lambda g.g\circ s)\phi \ast 0.\pi
_{X}$

Donc on se ram\`{e}ne \`{a} montrer :

$((\nu )\lambda g.g\circ s)\phi \ast 0.\pi _{X}\in \amalg $

Pour cela, on va interpr\'{e}ter judicieusement une variable de pr\'{e}dicat
$P$ qui servira \`{a} instancier $Int(n)$.

on d\'{e}finit $\Vert Pk\Vert =\{s^{n-k}0.\pi _{X}/\pi _{X}\in \Vert X\Vert
\}$ pour $0\leq k\leq n$ et $\varnothing $ sinon.

si on avait

a) $\lambda g.g\circ s\Vdash \forall x(Px\rightarrow Psx)$

b) $\phi \Vdash P0$

on aurait gagn\'{e}. En effet, $\nu \Vdash \forall x(Px\rightarrow
Psx),P0\rightarrow Pn,$ d'o\`{u} $\ $on tire $\nu \ast \lambda g.g\circ
s.\phi .0.\pi _{X}\in \amalg .$ Comme $((\nu )\lambda g.g\circ s)\phi \ast
0.\pi _{X}\succ \nu \ast \lambda g.g\circ s.\phi .0.\pi _{X},$ la
cl\^{o}ture de $\amalg $ par anti-r\'{e}duction nous donne le r\'{e}sultat.

a) ?$\lambda g.g\circ s\Vdash \forall x(Px\rightarrow Psx)$

? $\lambda g.g\circ s\ast t_{Pk}.\pi _{Psk}\in \amalg $

si $k+1\succ n,$c'est termin\'{e}, car alors $\left| Psk\right| =\Lambda ,$
et donc $\lambda g.g\circ s.t_{Pk}\in $ $\left| Psk\right| $

sinon $\pi _{Psk}\equiv s^{n-(k+1)}0.\pi _{X}$

? $\lambda g\lambda xg(s)x\ast t_{Pk}.s^{n-(k+1)}0.\pi _{X}\in \amalg $

? $\lambda xt_{Pk}(s)x\ast s^{n-(k+1)}0.\pi _{X}\in \amalg $

?$t_{Pk}.(s)s^{n-(k+1)}0.\pi _{X}\in \amalg $

Ce qui est vrai car $s^{n-k}0.\pi _{X}\in \Vert Pk\Vert .$

b) Ceci d\'{e}coule imm\'{e}diatement de l'hypoth\`{e}se. En effet, $\phi
\Vdash P0$ si et seulement si $\phi \ast s^{n}0.\pi _{X}\in \amalg .$

On montre maintenant le th\'{e}or\`{e}me \ref{5}. On se place dans le cas
d'une fonction r\'{e}cursive $f$\ unaire. On pose $\hat{k}=_{def}s^{k}0$ et
on utilise la notion de repr\'{e}sentation habituelle pour les fonctions
r\'{e}cursives dans le lambda-calcul : un $\lambda $-terme $\phi $
repr\'{e}sente une fonction r\'{e}cursive totale f si, pour tout n, f(n)=p
implique $\phi \hat{n}=_{\beta }\lambda f\lambda x(f)^{p}x$ $p=f(n).$

Soit $\phi $ un $\lambda $-terme repr\'{e}sentant $f$ et $n$ un entier. $%
\phi \hat{n}=_{\beta }\lambda f\lambda x(f)^{p}x$ $\ o\grave{u}$ $p=f(n)$
donc $\phi \hat{n}>\lambda f.\xi $ avec $\lambda f.\xi =_{\beta }\lambda
f\lambda x(f)^{p}x.$ Par le th\'{e}or\`{e}me \ref{8}, $\lambda f.\xi \Vdash
Int[s^{p}0].$ Soit $\pi _{Int[s^{p}0]}\in \left\| Int[s^{p}0]\right\| ,$ on
a alors $\lambda f.\xi \ast \pi _{Int[s^{p}0]}\in \amalg .$

Par le lemme \ref{6}, $\phi \hat{n}\ast \pi \succ \lambda f.\xi \ast \pi ,$
et comme $\phi \hat{n}$ est une application, ce que n'est pas $\lambda f.\xi
,$ la r\'{e}duction commence par $\phi \hat{n}\ast \pi \succ \phi \ast \hat{n%
}.\pi .$ Donc $\phi \ast \hat{n}.\pi \succ \lambda f.\xi \ast \pi $ et en
particulier $\phi \ast \hat{n}.\pi _{Int[s^{p}0]}\succ \lambda f.\xi \ast
\pi _{Int[s^{p}0]}\in \amalg .$ Ceci correspond \`{a} l'hypoth\`{e}se du
lemme \ref{9} pour $X=$.$Int[s^{p}0].$ Donc $T\phi \Vdash
Int[s^{n}0]\rightarrow Int[s^{p}0].$ Comme nous l'avons pour un $n$
quelconque, on a bien $T\phi \Vdash \forall x(Int[x]\rightarrow Int[f(x)]),$
autrement dit la formule $\forall x(Int[x]\rightarrow Int[f(x)])$ est
r\'{e}alis\'{e}e.

\subsection{La sp\'{e}cification des th\'{e}or\`{e}mes de l'arithm\'{e}tique}

\subsubsection{Pour les \'{e}nonc\'{e}s $\Pi _{2}$}

On montre que la preuve d'un \'{e}nonc\'{e} $\Pi _{2}$ de l'arithm\'{e}tique
dans la logique classique du second ordre peut \^{e}tre vu comme un
programme qui calcule la fonction associ\'{e}e \`{a} l'\'{e}nonc\'{e}.

\begin{theorem}
Si $\vdash \theta :\left[ \forall x\exists y(f(x,y)=0)\right] ^{int}$, alors
$\theta \ast \widehat{n}.Tt.\pi $ s'\'{e}value sur $t\ast \widehat{p}\pi
^{\prime }$ avec $f(n,p)=0$ o\`{u} $t\equiv \lambda x\lambda y.yx$
\end{theorem}

Le lemme d'ad\'{e}quation nous dit que $\theta \Vdash \forall
xInt(x)\rightarrow \left[ \exists y(f(x,y)=0)\right] ^{int}$. Donc pour un
entier $n$ quelconque,

$\theta \Vdash Int(n)\rightarrow \left[ \exists y(f(n,y)=0)\right] ^{int}$.
On fixe $\amalg =\{t\ast \widehat{p}.\pi ^{\prime }/f(n,p)=0,\pi ^{\prime
}\in \Pi \}^{\succ ^{-1}}.$

Le th\'{e}or\`{e}me \ref{8} nous assure que $\widehat{n}\Vdash Int(n).$
Reste donc \`{a} montrer que $Tt.\pi \in \left\| \left[ \exists y(f(n,y)=0)%
\right] ^{int}\right\| .$

$\left[ \exists y(f(n,y)=0)\right] ^{int}\equiv \forall y\left(
Int(y)\rightarrow (f(n,y)=0\rightarrow \perp )\right) \rightarrow \perp .$
Donc il faut avoir

$Tt\Vdash Int(m)\rightarrow (f(n,m)=0\rightarrow \perp )$ pour tout entier $%
m.$

Par le lemme \ref{9}, on se ram\`{e}ne \`{a} montrer

$t\ast \widehat{m}.\pi _{f(n,m)=0\rightarrow \perp }\in \amalg $ pour toute
pile $\pi _{f(n,m)=0\rightarrow \perp }\in \left\| f(n,m)=0\rightarrow \perp
\right\| $

Premier cas, $f(n,m)=0,$ dans ce cas, $t\ast \widehat{m}.\pi
_{f(n,m)=0\rightarrow \perp }\in \amalg $ par d\'{e}finition de $\amalg .$

Deuxi\`{e}me cas, $f(n,m)\neq 0.$ $\pi _{f(n,m)=0\rightarrow \perp }$ est de
la forme $u.\rho $ avec $u\in \left| \top \rightarrow \bot \right| $ et $%
\rho \in \Pi .$

$t\ast \widehat{m}.u.\rho \succ u\ast \widehat{m}.\rho $. Or $u\ast \widehat{%
m}.\rho \in \amalg $ car $\widehat{m}.\rho $ est trivialement dans $\left\|
\top \rightarrow \bot \right\| .$ Par cons\'{e}quent, on a bien $t\ast
\widehat{m}.\pi _{f(n,m)=0\rightarrow \perp }\in \amalg .$

La comparaison avec ce que donnent les m\'{e}thodes d'Ackermann et de
Kreisel est d\'{e}j\`{a} instructive pour ce cas.

Si $PA_{\epsilon }\vdash \forall x\exists yA(x,y),$ un contre-exemple serait
un $a$ tel que $\forall y\sim A(a,y).$ La NCI donne donc une fonction $f$
telle que pour tout $a,$ $A(a,f(a)),$ de sorte que $f$ et le code $\theta $
d'une preuve jouent bien le m\^{e}me r\^{o}le. N\'{e}anmoins, on remarque
que la fonction $f$ de Kreisel \'{e}crase l'individualit\'{e} des preuves :
le calcul de $f(a)$ commence par le calcul d'une limite sup\'{e}rieure pour
les valeurs de $f(a)$, cette limite \'{e}tant la m\^{e}me pour des preuves
de $\forall x\exists yA(x,y)$ qui ont les m\^{e}mes param\`{e}tres. En outre
la limite en question n'a pas de rapport direct avec la valeur finale de $%
f(a).$ Au contraire, $\theta $ est \`{a} chaque fois un programme singulier
qui calcule directement la valeur pour $f(a).$ En un sens, c'est d'ailleurs
bien exactement ce qu'est la m\'{e}thode de substitution : appliqu\'{e}e aux
formules de la preuve, elle peut-\^{e}tre vue comme un algorithme pour
calculer les valeurs de $f(a).$

\subsubsection{Pour les \'{e}nonc\'{e}s $\Sigma _{2}$}

On va montrer en quoi le terme codant la preuve d'un \'{e}nonc\'{e} $\Sigma
_{2}$ constitue un fonctionnel analogue aux fonctionnels r\'{e}cursifs de
Kreisel. La d\'{e}monstration suit celle de \cite{kr2}.

On dira qu'un terme $t$ repr\'{e}sente une fonction $f$ calculable \`{a} un
argument si

$t\ast \hat{n}.\pi \succ \widehat{f(n)}\ast \pi $

Le r\^{o}le de $t$ peut \^{e}tre jou\'{e} soit par un v\'{e}ritable $\lambda
$-terme, soit par une constante $c_{f}$ \'{e}quip\'{e}e de la r\`{e}gle de
r\'{e}duction ad\'{e}quate. Afin de r\'{e}soudre les difficult\'{e}s
pos\'{e}es par la d\'{e}finition de la r\'{e}duction dans le $\lambda _{c}$%
-calcul, on introduit une nouvelle constante $\zeta $ dont on d\'{e}finit
ainsi le comportement face \`{a} une pile.

Si $t\ast \hat{n}.\pi \succ \hat{p}\ast \pi ,$ alors pour tout $\xi ,$ $%
\zeta \ast \xi .(t)\hat{n}.\pi ^{\prime }\succ \xi \ast \hat{p}.\pi ^{\prime
}$

On pose $F[f]=(T)\lambda x\lambda y(((\zeta )y)(f)x)x$

\begin{theorem}
Si $\theta $ est le code d'une preuve de $[\exists x\forall y(\phi
(x,y)=0)]^{Int}$ dans l'arithm\'{e}tique du second ordre et $\gamma $ une
fonction r\'{e}cursive totale repr\'{e}sent\'{e}e par un terme t, alors,
pour toute pile $\pi $, le programme $\lambda f\theta F[f]$ auquel on donne
l'entr\'{e}e $t.\pi $ aboutit \`{a} un \'{e}tat $\hat{n}\ast \pi ^{\prime }$
tel que $N\vDash \phi (n,\gamma (n))=0$.
\end{theorem}

On commence par fixer $\amalg :\amalg =$\{ $\hat{n}\ast \pi ^{\prime }$ / $%
N\vDash \phi (n,\gamma (n))=0,\pi ^{\prime }\in \Pi $\}$^{\succ ^{-1}}$

Par le lemme d'ad\'{e}quation, il suffit de montrer la propri\'{e}t\'{e}
pour un $\theta $ qui r\'{e}alise $[\exists x\forall y(\phi (x,y)=0)]^{Int}.$
On a

$\theta \Vdash \forall x[Int(x),\forall y(Int(y)\rightarrow \phi
(x,y)=0)\rightarrow \bot ]\rightarrow \bot $

et $\lambda f\theta F[f]\ast t.\pi \succ $ $\theta \ast F[f:=t].\pi $

Donc, comme $\amalg $ est clos par anti-r\'{e}duction, il suffit d'avoir

$F[f:=t]\Vdash \forall x[Int(x),\forall y(Int(y)\rightarrow \phi
(x,y)=0)\rightarrow \bot ]$

On va montrer, pour un $n$ quelconque :

$F[f:=t]\Vdash Int(n),\forall y(Int(y)\rightarrow \phi (n,y)=0)\rightarrow
\bot $

Or les piles qui sont dans $\left\| \forall y(Int(y)\rightarrow \phi
(n,y)=0)\rightarrow \bot \right\| $ sont de la forme $\xi .\pi $ pour une
pile $\pi $ quelconque et un $\xi $ tel que $\xi \Vdash \forall
y(Int(y)\rightarrow \phi (n,y)=0).$ Donc par le lemme \ref{9}, on se
ram\`{e}ne \`{a} montrer~:

$F[f:=t]\ast \hat{n}.\xi .\pi \in \amalg $ c'est-\`{a}-dire $\lambda
x\lambda y(((\zeta )y)(t)x)x\ast \hat{n}.\xi .\pi \in \amalg $

Par ailleurs,

$\lambda x\lambda y(((\zeta )y)(t)x)x\ast \hat{n}.\xi .\pi $

$\succ (((\zeta )\xi )(t)\hat{n})\hat{n}\ast \pi $

$\succ \zeta \ast \xi .(t)\hat{n}.\hat{n}.\pi $

$\succ \xi \ast \hat{p}.\hat{n}.\pi $ o\`{u} $\gamma (n)=p,$ puisque $t\ast
\hat{n}.\pi \succ \widehat{f(n)}\ast \pi .$

Il suffit donc de montrer maintenant que $\xi \ast \hat{p}.\hat{n}.\pi $ est
dans $\amalg $. Comme $\xi \Vdash \forall y(Int(y)\rightarrow \phi (n,y)=0)$%
, $\xi .\hat{p}\Vdash \phi (n,p)=0,$ cela revient \`{a} \'{e}tablir que $%
\hat{n}.\pi \in \Vert \phi (n,p)=0\Vert $

Premier cas : $\phi (n,p)=0$.

On a alors $\Vert \phi (n,p)=0\Vert =\Vert \forall Z,Z\rightarrow Z\Vert $.

$\Vert \forall Z,Z\rightarrow Z\Vert \supseteq |\Pi |\rightarrow \Vert \Pi
\Vert $, or on a automatiquement $\pi \in $ $\Vert \Pi \Vert $ et, par
d\'{e}finition de $\amalg $, $\hat{p}\in |\Pi |$. Il suit que $\hat{n}.\pi
\in \Vert \phi (n,p)=0\Vert .$

Deuxi\`{e}me cas : $\phi (n,p)\neq 0$.

On a alors $||\phi (n,p)=0||=\top \rightarrow \bot $, ce qui nous donne
imm\'{e}diatement $\hat{n}.\pi \in \Vert \phi (n,p)=0\Vert $

CQFD

On constate que dans la d\'{e}monstration, la seule propri\'{e}t\'{e} de $%
F[t]$ qui est utilis\'{e}e est

$F^{\prime }[t]\ast \hat{n}.\xi \pi \succ \xi \ast \hat{p}.\hat{n}.\pi $
o\`{u} $F^{\prime }[t]$ est $F[t]$ moins le $T$ de l'application principale.

Autrement dit $F^{\prime }[t]$ doit fournir en r\'{e}ponse \`{a} un entier $%
n $ un autre entier $p$ qui correspond \`{a} la tentative d'exhiber un
contre-exemple et mettre en t\^{e}te le terme $\xi $ qui assure le
traitement par le programme de cette r\'{e}ponse. Intuitivement, $\xi $ va
ensuite mettre en t\^{e}te $\hat{n}$ si $n$ marche et proposer un autre $%
\hat{n}$ sinon. Mais il n'y a pas de raison de consid\'{e}rer que $F^{\prime
}[t]$ correspond \`{a} une strat\'{e}gie de r\'{e}futation d\'{e}cid\'{e}e
\`{a} l'avance correspondant \`{a} une fonction r\'{e}cursive. C'est dans
cet esprit qu'est formul\'{e} le th\'{e}or\`{e}me dans \cite{kr}\ et \cite
{kr2}. On ajoute aux termes du lambda-c calcul une constante $\kappa $ pour
laquelle on fixe un comportement ad\'{e}quat : $\kappa \ast \hat{n}.\xi .\pi
\gg \xi \ast \hat{p}.[\hat{n},\hat{p}].\pi $ ($[\hat{n},\hat{p}]$ se laisse
construire comme une constante ou comme une liste d'entiers); $\kappa $ peut
\^{e}tre vu comme une instruction \textit{input} demandant une r\'{e}ponse
\`{a} $n.$. Le th\'{e}or\`{e}me s'\'{e}nonce alors :

Si $\theta $ est le code d'une preuve de $[\exists x\forall y(\phi
(x,y)=0)]^{Int}$ dans l'arithm\'{e}tique du second ordre, alors, pour toute
pile $\pi $, le programme $\theta $ auquel on donne l'entr\'{e}e $T\kappa
.\pi $ aboutit \`{a} un \'{e}tat $[\hat{n},\hat{p}]\ast \pi ^{\prime }$ tel
que $N\vDash \phi (n,p)=0$.

Remarquons n\'{e}anmoins qu'au cours d'une r\'{e}duction, $\kappa $ n'arrive
qu'un nombre fini de fois en t\^{e}te, de sorte que les valeurs donn\'{e}es
par l'instruction \textit{input }peuvent toujours en droit \^{e}tre vues
comme produites par une fonction r\'{e}cursive.

On peut alors voir le programme comme une strat\'{e}gie gagnante pour $%
\exists loise$ dans un jeu P s\'{e}mantique suivant : le jeu commence par le
choix d'un entier $x$\ par $\exists loise,$ puis continue par un choix d'un
entier $y$ par $\forall b\acute{e}lard.$ $\exists loise$ peut alors soit
s'arr\^{e}ter et $\exists loise$ gagne si $\phi (x,y)=0$ tandis qu'$\forall b%
\acute{e}lard$ perd dans le cas contraire. Mais $\exists loise$ peut
\'{e}galement choisir de recommencer \`{a} proposer un $x.$ Si le jeu dure
infiniment longtemps, $\forall b\acute{e}lard$ gagne.

Une ex\'{e}cution du programme sur $T\kappa .\pi $ revient au fait de jouer
une partie, l'\'{e}tat $\kappa \ast \hat{n}.\xi .\pi $ correspond au choix d'%
$\forall b\acute{e}lard$ et tous les \'{e}tats $\xi \ast \hat{p}.\hat{n}.\pi
$ atteints lors de la r\'{e}duction correspondent au fait d'atteindre une
des positions finales.

La m\'{e}thode de substitution peut elle aussi \^{e}tre interpr\'{e}t\'{e}e
comme fournissant une strat\'{e}gie gagnant contre toute strat\'{e}gie de
l'adversaire. On commence par annuler tous les substituants des $\epsilon $%
-termes, ce qui revient \`{a} proposer 0; alors de deux choses l'une, soit
tous les axiomes de la forme III.1 sont vrais, et 0 est gagnant, soit ce
n'est pas le cas, l'algorithme continue alors \`{a} tourner en modifiant les
substituants des cat\'{e}gories pour rendre vrais les axiomes de la forme
III.1 - ce qui ne se traduit pas forc\'{e}ment \`{a} chaque pas par une
modification de la valeur finale propos\'{e}e. Ce que la preuve de
terminaison \'{e}tablit, c'est qu'au bout d'un moment, l'algorithme rend
vrai tous les axiomes critiques, ce qui implique que le substituant
propos\'{e} pour l'$\epsilon $-terme de la formule finale donne une valeur
gagnante$.$ On pourrait penser que si la strat\'{e}gie de l'adversaire n'est
pas identifi\'{e}e d\`{e}s le d\'{e}part \`{a} une fonction r\'{e}cursive,
on ne peut pas donner \`{a} l'avance une limite au nombre de substitutions,
mais en fait, on a vu dans le lemme \ref{12} que la fonctionnelle qui
fournit cette limite ne d\'{e}pend que d'un nombre fini de valeurs des
fonctions, ce qui veut dire qu'on n'a besoin que d'un nombre fini de tests
sur les r\'{e}ponses de l'adversaire pour borner la valeur. Par contre, rien
ne garantit que l'on s'arr\^{e}te d\`{e}s que l'on sort une valeur gagnante
(il se peut que des axiomes soient encore faux dans la preuve), il en va
d'ailleurs de m\^{e}me avec les programmes que constituent les $\lambda _{c}$%
-termes.

\subsubsection{Dans le cas g\'{e}n\'{e}ral}

On dira qu'un terme $t$ repr\'{e}sente une fonction $f$ calculable \`{a} k
arguments si

$t\ast \widehat{n_{1}}...\widehat{n_{k}}.\pi \succ \widehat{f(n_{1}...n_{k})}%
\ast \pi .$

Le r\^{o}le de $t$ peut \^{e}tre jou\'{e} soit par un v\'{e}ritable $\lambda
$-terme, soit une constante $c_{f}$ \'{e}quip\'{e}e de la r\`{e}gle de
r\'{e}duction ad\'{e}quate. De mani\`{e}re analogue au $\zeta $
pr\'{e}c\'{e}dent, on d\'{e}finit des $\zeta _{k}$ qui permettent
d'\'{e}valuer des fonctions \`{a} k arguments \`{a} l'int\'{e}rieur d'une
pile.

Si $t\ast \widehat{n_{1}}...\widehat{n_{k}}.\pi \succ \widehat{p}\ast \pi ,$
alors pour tout $\xi ,$ $\zeta _{k}\ast \xi .(t)\widehat{n_{1}}...\widehat{%
n_{k}}.\pi ^{\prime }\succ \xi \ast \hat{p}.\pi ^{\prime }$

On commence par d\'{e}finir une suite de termes par induction r\'{e}trograde
:

$H_{k}=[x_{1},...,x_{k}]$ pour une repr\'{e}sentation fix\'{e}e des listes
et pour un k donn\'{e}

$H_{j}=(T)\lambda x_{j+1}\lambda y_{j+1}.(((\zeta
_{j+1})y_{j+1})(f_{j+1})x_{1}...x_{j+1})H_{j+1}$ pour $0\leq j<k$

On d\'{e}signe alors par $F_{k}[f_{1}...f_{k}]$ le terme $H_{0}.$ On
remarque que le $F$ de la section pr\'{e}c\'{e}dente correspond bien \`{a}
notre $H_{1}$ actuel.

\begin{theorem}
Si $\theta $ est le code d'une preuve de

$[\exists x_{1}\forall y_{1}...\exists x_{k}\forall y_{k}(\phi
(x_{1}...x_{k},y_{1}...y_{k})=0)]^{Int}$ dans l'arithm\'{e}tique du second
ordre et $\gamma _{1}(x_{1})...\gamma _{k}(x_{1}...x_{k})$ sont des
fonctions r\'{e}cursives totales repr\'{e}sent\'{e}es par les $\lambda $%
-terme $t_{1}...t_{k}$, alors, pour toute pile $\pi $, $\lambda
f_{1}...\lambda f_{k}\theta F_{k}[f_{1}...f_{k}]\ast t_{1}...t_{k}.\pi $ se
r\'{e}duit sur $[\hat{n}_{1},...,\hat{n}_{k}]\ast \pi ^{\prime }$ avec $%
N\vDash \phi (n_{1}...n_{k},\gamma _{1}(n_{1})...y_{k}(n_{1}...n_{k}))=0)$.
\end{theorem}

On commence par fixer $\amalg :\amalg =$\{ $[\hat{n}_{1},...,\hat{n}%
_{k}]\ast \pi $ / $N\vDash \phi (n_{1}...n_{k},\gamma
_{1}(n_{1})...y_{k}(n_{1}...n_{k}))=0$\}

Par le lemme d'ad\'{e}quation, il suffit de montrer la propri\'{e}t\'{e}
pour un $\theta $ qui r\'{e}alise $[\exists x_{1}\forall y_{1}...\exists
x_{k}\forall y_{k}\phi (x_{1}...x_{k},y_{1}...y_{k})=0]^{Int}.$

On pose

$A_{i}(x_{1}...x_{i})\equiv \lbrack \exists x_{i+1}\forall y_{i+1}...\exists
x_{k}\forall y_{k}\phi (x_{1}...x_{k},\gamma _{1}(x_{1})...\gamma
_{i}(x_{1}...x_{i})y_{i+1}...y_{k})=0]^{Int},$ $A_{i}^{\prime
}(x_{1}...x_{i})\equiv \lbrack \exists x_{i+1}\forall y_{i+1}...\exists
x_{k}\forall y_{k}\phi (x_{1}...x_{k},\gamma _{1}(x_{1})...\gamma
_{i-1}(x_{1}...x_{i-1})y_{i}y_{i+1}...y_{k})=0]^{Int}$, $H_{j}^{\prime
}(n_{1}...n_{j})\equiv H_{j}[x_{1}:=\hat{n}_{1},...,x_{j}:=\hat{n}%
_{j},f_{1}:=t_{1},...,f_{k}:=t_{k}].$

$\lambda f_{1}...\lambda f_{k}\theta F_{k}[f_{1}...f_{k}]\ast
t_{1}...t_{k}.\pi \succ \theta H_{0}^{\prime }\ast \pi \succ \theta \ast
H_{0}^{\prime }.\pi $ donc comme $\theta \Vdash \lbrack \exists x_{1}\forall
y_{1}...\exists x_{k}\forall y_{k}\phi
(x_{1}...x_{k},y_{1}...y_{k})=0]^{Int} $, il s'agit de montrer $%
H_{0}^{\prime }.\pi \in \Vert A_{0}\Vert .$ La d\'{e}monstration se fait par
induction r\'{e}trograde en commen\c{c}ant par $H_{k}^{\prime
}(n_{1}...n_{k}).\pi \in \Vert \phi (n_{1}...n_{k},\gamma
_{1}(n_{1})...y_{k}(n_{1}...n_{k}))=0)\Vert .$ La situation est semblable
\`{a} la fin de la d\'{e}monstration pr\'{e}c\'{e}dente.

Premier cas : $\phi (n_{1}...n_{k},\gamma
_{1}(n_{1})...y_{k}(n_{1}...n_{k}))=0$.

On a alors $\Vert \phi (n_{1}...n_{k},\gamma
_{1}(n_{1})...y_{k}(n_{1}...n_{k}))=0\Vert =\Vert \forall Z,Z\rightarrow
Z\Vert $.

$\Vert \forall Z,Z\rightarrow Z\Vert \supseteq |\Pi |\rightarrow \Vert \Pi
\Vert $, or on a automatiquement $\pi \in $ $\Vert \Pi \Vert $ et, par
d\'{e}finition de $\amalg $, $[\hat{n}_{1},...,\hat{n}_{k}]\in |\Pi |$. Il
suit que $[\hat{n}_{1},...,\hat{n}_{k}].\pi \in \Vert \lbrack \hat{n}%
_{1},...,\hat{n}_{k}]\Vert .$

Deuxi\`{e}me cas : $\phi (n_{1}...n_{k},\gamma
_{1}(n_{1})...y_{k}(n_{1}...n_{k}))\neq 0$

On a alors $||\phi (n,p)=0||=\top \rightarrow \bot $, ce qui nous donne
imm\'{e}diatement $[\hat{n}_{1},...,\hat{n}_{k}].\pi \in \Vert \phi
(n_{1}...n_{k},\gamma _{1}(n_{1})...y_{k}(n_{1}...n_{k}))=0\Vert $

On montre maintenant $H_{j}^{\prime }(n_{1}...n_{j}).\pi \in \Vert
A_{j}(n_{1}...n_{j})\Vert $ sous l'hypoth\`{e}se $H_{j+1}^{\prime
}(n_{1}...n_{j+1}).\pi \in \Vert A_{j+1}(n_{1}...n_{j+1})\Vert .$

On commence par remarquer que, pour un terme $\xi $ et une pile $\pi $
quelconques, on a :

$\lambda x_{j+1}\lambda y_{j+1}.(((\zeta _{j+1})y_{j+1})(t_{j+1})\hat{n}%
_{1}...\hat{n}_{j}x_{j+1})H_{j+1}^{\prime }(n_{1}...n_{j}x_{j+1})\ast \hat{n}%
_{j+1}.\xi .\pi $

$\succ (((\zeta _{j+1})\xi )(t_{j+1})\hat{n}_{1}...\hat{n}_{j}\hat{n}%
_{j+1})H_{j+1}^{\prime }(n_{1}...n_{j}\hat{n}_{j+1})\ast .\pi $

$\succ \zeta _{j+1}\ast \xi .(t_{j+1})\hat{n}_{1}...\hat{n}_{j}\hat{n}%
_{j+1}.H_{j+1}^{\prime }(n_{1}...n_{j}\hat{n}_{j+1}).\pi $ puis comme $%
t_{j+1}\ast \hat{n}_{1}...\hat{n}_{j+1}.\pi \succ \hat{p}_{j+1}.\pi $ o\`{u}
$p_{j+1}=\gamma (n_{1}...n_{j+1})$

$\succ \xi \ast \hat{p}_{j+1}.H_{j+1}^{\prime }(n_{1}...n_{j+1}).\pi $ \ \ \
\ (*)

On doit montrer $H_{j}^{\prime }(n_{1}...n_{j}).\pi \in \Vert \forall
x_{j+1}[Int(x_{j+1}),\forall y_{j+1}(Int(y_{j+1})\rightarrow
A_{j+1}(n_{1}...n_{j}x_{j+1}))\rightarrow \bot ]\rightarrow \perp \Vert ,$
ce qui revient \`{a} $H_{j}^{\prime }(n_{1}...n_{j})\Vdash \forall
x_{j+1}[Int(x_{j+1}),\forall y_{j+1}(Int(y_{j+1})\rightarrow
A_{j+1}(n_{1}...n_{j}x_{j+1}))\rightarrow \bot ]$ ou encore $H_{j}^{\prime
}(n_{1}...n_{j})\Vdash Int(n_{j+1})\rightarrow \lbrack \forall
y_{j+1}(Int(y_{j+1})\rightarrow A_{j+1}(n_{1}...n_{j+1}))\rightarrow \bot ]$
pour un $n_{j+1}$ quelconque.

$H_{j}^{\prime }(n_{1}...n_{j})\equiv (T)\lambda x_{j+1}\lambda
y_{j+1}.(((\zeta _{j+1})y_{j+1})(t_{j+1})\hat{n}_{1}...\hat{n}%
_{j}x_{j+1})H_{j+1}^{\prime }(n_{1}...n_{j}x_{j+1})$. Or les piles qui sont
dans \ $\left\| \forall y_{j+1}(Int(y_{j+1})\rightarrow
A_{j+1}(n_{1}...n_{j+1}))\rightarrow \bot \right\| $ sont de la forme $\xi
.\pi $ o\`{u} $\pi $ une pile quelconque et $\xi $ un terme r\'{e}alisant $%
\forall y_{j+1}(Int(y_{j+1})\rightarrow A_{j+1}^{\prime }(n_{1}...n_{j+1})).$
Donc, d'apr\`{e}s le lemme \ref{9} il suffit de montrer : $\lambda
x_{j+1}\lambda y_{j+1}.(((\zeta _{j+1})y_{j+1})(t_{j+1})\hat{n}_{1}...\hat{n}%
_{j}x_{j+1})H_{j+1}^{\prime }(n_{1}...n_{j}x_{j+1})\ast \hat{n}_{j+1}.\xi
.\pi \in \amalg $

Par (*) et la cl\^{o}ture de $\amalg $ par anti-r\'{e}duction, il suffit
pour cela d'avoir $\xi \ast \hat{p}_{j+1}.H_{j+1}^{\prime
}(n_{1}...n_{j+1}).\pi \in \amalg .$ Mais, par d\'{e}finition, $\xi $
r\'{e}alise $Int(p_{j+1})\rightarrow A_{j+1}(n_{1}...n_{j+1}),$ donc il
suffit de v\'{e}rifier que $\hat{p}_{j+1}\in |Int(n_{j+1})|$ et $%
H_{j+1}^{\prime }(n_{1}...n_{j+1}).\pi \in \left\|
A_{j+1}(n_{1}...n_{j+1})\right\| ,$ ce qui est imm\'{e}diat pour $\hat{p}%
_{j+1}$ et d\'{e}coule de l'hypoth\`{e}se d'induction pour $H_{j+1}^{\prime
}(n_{1}...n_{j+1}).\pi .$

CQFD.

On peut faire les m\^{e}mes remarques que pour le th\'{e}or\`{e}me
pr\'{e}c\'{e}dent. Ici, la propri\'{e}t\'{e} utile du terme universel $%
F_{k}[t_{1}...t_{k}]$ et de ses sous termes est (*) de sorte que l'on peut
faire jouer le r\^{o}le des $H_{j}$ par des constantes munies d'une
r\`{e}gle de r\'{e}duction appropri\'{e}es; \cite{kr} dote ainsi des
constantes $\kappa _{n_{1}p_{1},...,n_{j}p_{j}}^{j}$ des r\`{e}gles de
r\'{e}duction suivantes

pour $0\leq j\leq k-2,$ $\kappa _{n_{1}p_{1},...,n_{j}p_{j}}^{j}\widehat{n}%
_{j+1}\ast \xi .\pi \gg \xi \ast \widehat{p}_{j+1}.T\kappa
_{n_{1}p_{1},...,n_{j+1}p_{j+1}}^{j}$

pour $j=k-1,$ $\kappa _{n_{1}p_{1},...,n_{k-1}p_{k-1}}^{j}\widehat{n}%
_{k}\ast \xi .\pi \gg \xi \ast \widehat{p}_{k}.[\widehat{n}_{1}\widehat{p}%
_{1},...,\widehat{n}_{k}\widehat{p}_{k}]$

et montre le th\'{e}or\`{e}me :

Si $\theta $ est le code d'une preuve de

$[\exists x_{1}\forall y_{1}...\exists x_{k}\forall y_{k}(\phi
(x_{1}...x_{k},y_{1}...y_{k})=0)]^{Int}$ dans l'arithm\'{e}tique du second
ordre alors, pour toute pile $\pi $, $\theta \ast T\kappa ^{0}.\pi $ se
r\'{e}duit sur $[\widehat{n}_{1}\widehat{p}_{1},...,\widehat{n}_{k}\widehat{p%
}_{k}]\ast \pi ^{\prime }$ avec $N\vDash \phi
(n_{1}...n_{k},p_{1}...p_{k})=0. $

L\`{a} encore l'interpr\'{e}tation en termes de jeu est naturelle; le jeu
est le m\^{e}me que pr\'{e}c\'{e}demment; les quantificateurs existentiels
correspondent aux coups d'$\exists lo\ddot{\imath}se,$ les universels \`{a}
ceux d'$\forall b\acute{e}lard$ et $\exists lo\ddot{\imath}se$ peut revenir
en arri\`{e}re \`{a} tout moment. Les constantes $\kappa
_{n_{1}p_{1},...,n_{j}p_{j}}^{j}$ correspondent toujours \`{a} des
instructions \textit{input }qui demandent \`{a} $\forall b\acute{e}lard$ ce
qu'il veut jouer en fonction de la position qui a \'{e}t\'{e} atteinte dans
le jeu. Le d\'{e}roulement du programme $\theta $ auquel on donne comme
entr\'{e}e $T\kappa ^{0}.\pi $ o\`{u} $\pi $ ne comporte pas de constantes $%
\kappa $, correspond bien \`{a} un d\'{e}roulement du jeu selon les
r\`{e}gles. En effet, si une constante $\kappa
_{n_{1}p_{1},...,n_{j+1}p_{j+1}}^{j+1}$ arrive en t\^{e}te, c'est
n\'{e}cessairement que des constantes $\kappa ^{0}...\kappa
_{n_{1}p_{1},...,n_{j}p_{j}}^{j}$ sont pr\'{e}c\'{e}demment arriv\'{e}s en
t\^{e}te, respectant la r\`{e}gle selon laquelle $\exists lo\ddot{\imath}se$
ne peut repartir que de positions qui ont d\'{e}j\`{a} \'{e}t\'{e}
atteintes; en outre, cela montre bien que tous les coups successifs du jeu
sont jou\'{e}s, au sens o\`{u} avant de terminer sur $[\widehat{n}_{1}%
\widehat{p}_{1},...,\widehat{n}_{k}\widehat{p}_{k}],$ le programme a
n\'{e}cessairement mis en t\^{e}te au moins une fois chacun des $\kappa ^{i}.
$

\subsection{Remarque finale sur les deux approches}

L'analyse de Kreisel mettait l'accent sur la modularit\'{e} de la notion
d'interpr\'{e}tation. Le fait de v\'{e}rifier imm\'{e}diatement et de
mani\`{e}re constructive les clauses $\gamma )$ et $\delta )$ est m\^{e}me
ce qui constitue la sp\'{e}cificit\'{e} de la d\'{e}monstration de la NCI
fond\'{e}e sur l'$\epsilon $-calcul par rapport aux autres approches,
qu'elles partent de la preuve de consistance Gentzen ou d'une
interpr\'{e}tation fonctionnelle (voir par exemple \cite{avfef}). Kohlenbach
\cite{kohl} a montr\'{e} les clauses $\gamma )$ et $\delta )$ gr\^{a}ce au
couple traduction de G\"{o}del et interpr\'{e}tation fonctionnelle, mais au
prix de lourdes manipulations, puisqu'il s'agit de montrer que ces clauses
sont des th\'{e}or\`{e}mes dans des extensions de PA$^{\omega }$ et de
fournir ensuite une interpr\'{e}tation fonctionnelle pour ces
th\'{e}or\`{e}mes.

Il est naturel de poser les m\^{e}mes questions dans le cadre de
l'interpr\'{e}tation par des $\lambda _{c}$-termes.\ La r\'{e}ponse sera
diff\'{e}rente selon le niveau auquel on se place.

Si par interpr\'{e}tation d'une formule, on entend seulement un $\lambda
_{c} $-terme qui est le code d'une preuve, la modularit\'{e} de
l'interpr\'{e}tation est triviale. Pour l'analogue de la condition $\delta ),
$ etant donn\'{e}s deux termes $t$ et $u$ correspondant respectivement \`{a}
une preuve\ de $A$ et \`{a} une preuve de $A\rightarrow B,$ $(u)t$ donne
bien s\^{u}r une interpr\'{e}tation pour $B,$ moyennnant \'{e}ventuellement
les transformations li\'{e}es \`{a} la mise sous forme pr\'{e}nexe. De
m\^{e}me pour la condition $\gamma )$, \'{e}tant donn\'{e} un terme $u$
codant une preuve de $\sim A,$ on montre qu'il n'y a pas de termes $t$
codant une preuve de $A$. $(u)t$ serait alors typable avec $\perp .$ Par le
lemme d'ad\'{e}quation, on se ram\`{e}ne \`{a} prouver qu'il n'y a pas de $%
\lambda _{c} $-terme $v$ tel que $v\Vdash \forall XX$ pour tout choix de $%
\amalg .$ On consid\`{e}re en effet $\amalg =\{c\ast \rho \}^{\succ ^{-1}}$
et $\amalg ^{\prime }=\{c^{\prime }\ast \rho \}^{\succ ^{-1}}$ o\`{u} $c$ et
$c^{\prime }$ sont deux constantes fix\'{e}es diff\'{e}rentes et $\rho $ un
fonds de pile. Soit $\pi \ $une pile quelconque, on devrait avoir \`{a} la
fois $v\ast \pi \succ c\ast \rho $ et $v\ast \pi \succ c^{\prime }\ast \rho $
ce qui est impossible. Le fait que $\perp $ ne soit pas r\'{e}alis\'{e}
constitue une preuve que $\perp $ n'est pas d\'{e}rivable dans AP$_{2},$ de
sorte qu'on retrouve ce dont la NCI partait, \`{a} savoir une preuve de
consistanc (m\^{e}me si les moyens utilis\'{e}s ici n'ont rien de finitistes
bien s\^{u}r).

Reste qu'en un sens, un $\lambda _{c}$-terme qui se comporte comme le code
d'une preuve, constitue tout aussi bien une interpr\'{e}tation que le code
d'une preuve. Soit $A\equiv \lbrack \exists x_{1}\forall y_{1}...\exists
x_{k}\forall y_{k}(\phi (x_{1}...x_{k},y_{1}...y_{k})=0)]^{Int}$ et $t$ un $%
\lambda _{c}$-terme tel que $t\Vdash A$ avec $\amalg _{A}=$\{$[\widehat{n}%
_{1}\widehat{p}_{1},...,\widehat{n}_{k}\widehat{p}_{k}]\ast \pi $ / $N\vDash
\phi (n_{1}...n_{k},\gamma _{1}(n_{1})...y_{k}(n_{1}...n_{k}))=0$\}. $t$
fournit une interpr\'{e}tation pour A, c'est-\`{a}-dire une strat\'{e}gie
gagnante dans le jeu associ\'{e} \`{a} A. La question pos\'{e}e par $\delta
) $ est alors : \'{e}tant donn\'{e}s une preuve de $A\rightarrow B$ et un $%
\lambda _{c}$-terme qui r\'{e}alise $A$ au sens de $\amalg _{A},$ peut-on
trouver un terme qui r\'{e}alise $B$ au sens de $\amalg _{A}$ ? On pourrait
encore g\'{e}n\'{e}raliser et demander : \'{e}tant donn\'{e}s deux $\lambda
_{c}$-termes qui r\'{e}alisent respectivement $A$ au sens de $\amalg _{A}$
et $A\rightarrow B$ au sens de $\amalg _{A\rightarrow B}$ peut-on trouver un
terme qui r\'{e}alise $B$ au sens de $\amalg _{A}$ ?

L'id\'{e}e est que l'interpr\'{e}tation de $A\rightarrow B$ doit \^{e}tre
une strat\'{e}gie gagnante dans le jeu qui commence par le choix par $%
\exists lo\ddot{\imath}se$ de $\sim A$ ou de $B$ et qui se poursuite ensuite
comme les jeux associ\'{e}s aux formes pr\'{e}nexes de $\sim A$ et $B.$ Il
semble alors n\'{e}cessaire de fournir une sp\'{e}cification pour $%
A\rightarrow B,$ qui n'est pas en forme pr\'{e}nexe; en fait on n'est pas
oblig\'{e} de r\'{e}soudre le probl\`{e}me directement. On peut se contenter
d'interpr\'{e}ter les connecteurs comme des quantificateurs. On peut
supposer sans perte de g\'{e}n\'{e}ralit\'{e} que $\sim A$ et $B$ sous forme
pr\'{e}nexe sont caract\'{e}ris\'{e}s par les m\^{e}mes alternance de
quantificateurs (au besoin, on ajoute des quantificateurs qui ne lient
aucune variable). Soient $\sim A\equiv \lbrack \forall x_{0}\exists
x_{1}\forall y_{1}...\exists x_{k}\forall y_{k}(\phi
(x_{0}...x_{k},y_{1}...y_{k})=0)$ et $B\equiv \lbrack \forall x_{0}\exists
y_{1}\forall y_{1}...\forall x_{k}\exists y_{k}(\psi
(x_{0}...x_{k},y_{1}...y_{k})=0)],$ une preuve de $A\rightarrow B$ nous
donne une preuve de $A\Rightarrow B\equiv \exists x_{0}\forall
x_{0}...\forall x_{k}\exists y_{k}(\vartheta
(x_{0}...x_{k},y_{0}...y_{k})=0) $ avec

$\vartheta (x_{0}...x_{k},y_{0}...y_{k})=0$ si et seulement si

$y_{0}=0$ et $\phi (x_{0}...x_{k},y_{1}...y_{k})=0$

$y_{0}\neq 0$ et $\psi (x_{0}...x_{k},y_{1}...y_{k})=0$

Le th\'{e}or\`{e}me de sp\'{e}cification nous dit bien alors que\ les
r\'{e}alisateurs de $A\Rightarrow B$ donnent une strat\'{e}gie gagnante dans
le jeu associ\'{e} \`{a} $A\rightarrow B.\cite{kr}$

Pour revenir \`{a} la question pos\'{e}e, Coquand \cite{coq}, dans lequel
l'\'{e}limination des coupures est vue comme une interaction entre deux
strat\'{e}gies sur des jeux tout-\`{a}-fait analogues aux jeux ici
consid\'{e}r\'{e}s,\ fournit une r\'{e}ponse partielle. Un algorithme est en
effet d\'{e}fini qui, \'{e}tant donn\'{e}es une formule $B$ existentielle et
des strat\'{e}gies gagnantes $\sigma $ et $\tau $ pour les formules $%
A\rightarrow B$ et $A,$ fournit un t\'{e}moin pour $B.$ L'id\'{e}e est
d'organiser un \textit{d\'{e}bat }qui fait jouer la strat\'{e}gie $\sigma $
vue comme strat\'{e}gie partielle sur $\sim A$ contre $\tau $ et de montrer
que le d\'{e}bat termine, c'est-\`{a}-dire qu'au bout d'un moment $\sigma $
est oblig\'{e} d'abandonner $\sim A$ et de jouer sur $B.$

\end{document}